\documentclass{amsart}

\usepackage{a4wide,amssymb,amsrefs,calligra,comment,graphicx,mathrsfs,txfonts}

\makeatletter

\let\c@equation=\c@subsection

\def\thmhead#1#2#3{%
  (\thmnumber{%\@ifnotempty{#1}{ }
    \@upn{#2}})
   \thmname{#1}%
   \thmnote{ {\the\thm@notefont(#3)}}}

\makeatother

\def\subsection#1{\refstepcounter{subsection}
\medskip\noindent{\textbf{(\thesubsection)\ #1\unskip. }}\ignorespaces}

\def\subsubsection#1{\refstepcounter{subsubsection}\smallskip\noindent
  {\textbf{(\thesubsubsection)}\ \textit{#1\unskip. }}\ignorespaces}

\def\thesubsection
 {\thesection.\arabic{subsection}}
\def\thesubsubsection
 {\thesection.\arabic{subsection}.\arabic{subsubsection}}

\def\theodef#1{\newtheorem{#1}[subsection]{#1}}

\usepackage[2cell,arrow,curve,matrix]{xy}
\UseHalfTwocells
\UseTwocells

\def\calli#1{\textup{\!\textcalligra{#1}\,}}

\def\xymat#1{\begin{aligned}\xymatrixcolsep{3pc}\xymatrix{#1}\end{aligned}}

\def\A{{\mathbf A}}
\def\a{{\mathbf a}}
\def\B{{\mathbf B}}
\def\b{{\mathbf b}}

\def\F{{\mathbb F}}

\def\g{{\mathbf g}}
\def\i{{\bar i}}
\def\L{{\mathrm L}}
\def\LL{{\mathbb L}}
\def\p{{\bar p}}
\def\R{{\mathrm R}}
\def\RR{{\mathbb R}}
\def\r{{\mathscr R}}
\def\t{{\mathscr T}}
\def\Z{{\mathbb Z}}

\def\op{^{\mathrm{op}}}
\def\ho{_\simeq}
\def\1{^{-1}}
\def\comp{{\scriptstyle\Box}}
\def\binv{^{\scriptscriptstyle\boxminus}}
\def\Mod{{\textbf{-Mod}}}
\def\mod{{\textbf{-mod}}}

\DeclareMathOperator{\Aut}{Aut}
\DeclareMathOperator{\Ext}{Ext}
\DeclareMathOperator{\Hom}{Hom}
\DeclareMathOperator{\im}{im}

\def\xto#1{\xrightarrow[]{#1}}
\def\xot#1{\xleftarrow[]{#1}}
\def\set#1{\left\{#1\right\}}
\def\setof#1#2{\left\{#1\ \mid\ #2\right\}}
\def\hog#1{\left\llbracket#1\right\rrbracket}

\def\brk#1{\left\langle#1\right\rangle}
\def\alignbox#1{\begin{aligned}#1\end{aligned}}

\let\d\delta
\let\ge\geqslant
\let\le\leqslant
\let\x\times
\let\ox\otimes
\let\ot\leftarrow
\let\then\Rightarrow

\let\toto\rightrightarrows

\theodef{Lemma}
\theodef{Proposition}
\theodef{Theorem}

\theoremstyle{definition}

\theodef{Definition}

\theodef{Remark}
\theodef{Example}
\theodef{Examples}

\begin{document}

\title
[Secondary derived functors and the Adams spectral sequence]
{Secondary derived functors\\and the Adams spectral sequence}

\author{Hans Joachim Baues}
\address{
Max-Planck-Institut f\"ur Mathematik\\
Vivatsgasse 7\\
D-53111 Bonn\\
Germany}
\email{baues@mpim-bonn.mpg.de}

\author{Mamuka Jibladze}
\address{
Razmadze Mathematical Institute\\
M.~Alexidze st. 1\\
Tbilisi 0193\\
Georgia}
\curraddr{
Max-Planck-Institut f\"ur Mathematik\\
Vivatsgasse 7\\
D-53111 Bonn\\
Germany}
\email{jib@rmi.acnet.ge}

\maketitle

The paper introduces the secondary derived functor $\Ext^n$ obtained by
secondary resolutions. This generalizes the concept of the classical
derived functor $\Ext^n$. The classical groups $\Ext^n$, for example, are
used to describe the E$_2$-term of the Adams spectral sequence. As a main
application we show that the secondary $\Ext$ groups determine the
E$_3$-term of the Adams spectral sequence. Using the theory in
\cite{Baues} this yields an algorithm for the computation of the
E$_3$-term, see \cite{E3}. The algorithm is achieved by taking into account the
\emph{track structure}: one considers not just homotopy classes of maps
between spectra, but instead maps and homotopy classes of homotopies
between maps, termed tracks. These form a \emph{track category}, that is,
a category enriched in groupoids. It then turns out that in appropriate
track categories secondary $\Ext$ groups can be defined which are
unchanged if one replaces the ambient track category with a weakly
equivalent one. In fact in \cite{Baues} a manageable purely algebraically
described track category weakly equivalent to the track category of
Eilenberg-Mac Lane spectra has been completely determined. It is this
algebraic model that will be used on the basis of the main result
\ref{coinc} below to compute explicitly the E$_3$-term
of the Adams spectral sequence as a secondary $\Ext$-group, see \cite{E3}.

\section{Derived functors}
\label{primary}

We first recall the notion of a resolution in an additive category from
which we deduce (primary) derived functors. Later we introduce the
secondary version of these notions in the context of an ``additive track
category'', see section \ref{atrack}.

Our initial data consist of an additive category $\A$ and a full additive
subcategory $\a$ of $\A$. The basic situation to have in mind is the
category $R\Mod$ of modules over a ring $R$ and its subcategory $R\mod$ of
free (or projective) $R$-modules. As another motivating example, coming
from topology, one considers for $\A$ the opposite of the stable homotopy
category and for $\a$ its full subcategory on objects represented by
finite products of Eilenberg-Mac Lane spectra over a fixed prime field
$\F_p$; then $\a$ is equivalent to the category of finitely generated free
modules over the mod $p$ Steenrod algebra.

\begin{Definition}\label{1-complex}
A \emph{chain complex} $(A,d)$ in $\A$ is a sequence of objects and morphisms
$$
...\to A_{n+1}\xto{d_n}A_n\xto{d_{n-1}}A_{n-1}\to...
$$
from $\A$, with $d_{n-1}d_n=0$ ($n\in\Z$).

A \emph{chain map} $f:(A,d)\to(A',d')$ is a sequence of morphisms
$f_n:A_n\to A'_n$ with $f_nd_n=d'_nf_{n+1}$, $n\in\Z$. For two maps
$f,f':(A,d)\to(A',d')$, a \emph{chain homotopy} $h$ from $f$ to $f'$ is a
sequence of morphisms $h_n:A_{n-1}\to A'_n$ satisfying
$f'_n=f_n+d'_nh_{n+1}+h_nd_{n-1}$, $n\in\Z$.

A chain complex $(A,d)$ is called \emph{$\a$-exact} if for any object $X$ from
the subcategory $\a$ the (ordinary) chain complex $\Hom_\A(X,A_\bullet)$ of abelian groups
$$
...\to\Hom_\A(X,A_{n+1})\xto{\Hom_\A(X,d_n)}\Hom_\A(X,A_n)\xto{\Hom_\A(X,d_{n-1})}\Hom_\A(X,A_{n-1})\to...
$$
is acyclic, i.~e. an exact sequence. Explicitly, this means that for any
$n\in\Z$, any object $X$ from $\a$ and any morphism $a_n:X\to A_n$ with
$d_{n-1}a_n=0$ there exists a morphism $a_{n+1}:X\to A_{n+1}$ with
$a_n=d_na_{n+1}$.

A chain map $f:A\to A'$ is an \emph{$\a$-equivalence} if for every $X$ in
$\a$ the chain map $\Hom_\A(X,f)$ is a quasiisomorphism. Thus a chain
complex $(A,d)$ is $\a$-exact if and only if the map $(0,0)\to(A,d)$ is an
$\a$-equivalence.
\end{Definition}

\begin{Definition}\label{1-resolution}
For an object $A$ of $\A$, an \emph{$A$-augmented chain complex}
$A_\bullet^\epsilon$ is a chain complex of the form
$$
...\to A_1\to A_0\to A\to0\to0\to...,
$$
i.~e. with $A_{-1}=A$ and $A_{-n}=0$ for $n>1$. We will consider such an
augmented chain complex as a map between chain complexes,
$\epsilon:A_\bullet\to A$, where $A_\bullet$ is the complex $...\to A_1\to
A_0\to0\to0\to...$ whereas $A$ is considered as a complex concentrated in
degree 0, with $\epsilon=d_{-1}:A_0\to A$ called the augmentation.

An \emph{$\a$-resolution} of $A$ is an $\a$-exact $A$-augmented chain
complex such that all $A_n$ for $n\ge0$ belong to $\a$. Thus an
$\a$-resolution $A_\bullet^\epsilon$ of an object $A$ is the same as a
chain complex $A_\bullet$ in $\a$ together with an $\a$-equivalence
$\epsilon:A_\bullet\to A$.

There are obvious dual notions of an $A$-coaugmented complex and
$\a$-coresolution of $A$. Namely, this means a complex (resp. $\a$-exact
complex) with $A_1=A$ and $A_n=0$ for $n>1$.
\end{Definition}

\begin{Lemma}\label{1-existence}
Suppose
\begin{itemize}
\item
coproducts of any families of objects of $\a$ exist in $\A$ and belong to
$\a$ again;
\item
there is a small subcategory $\g$ of $\a$ such that every object of $\a$
is a retract of a coproduct of a family of objects from $\g$.
\end{itemize}
Then every object of $\A$ has an $\a$-resolution.
\end{Lemma}

\begin{proof}
We begin by taking
$$
A_0=\coprod_{\substack{G\in\g\\a:G\to A}}G,
$$
with the obvious map $\epsilon:A_0\to A$ having $a$ for the $a$-th
component. Next we take
$$
A_1=\coprod_{\substack{G\in\g\\t_0:G\to A_0\\\epsilon t_0=0}}G,
$$
with a similar map $d_0:A_1\to A_0$ whose $t_0$-th component is $t_0$ (so
obviously $\epsilon d_0=0$). One continues in this way, with
$$
A_{n+1}=\coprod_{\substack{G\in\g\\t_n:G\to A_n\\d_{n-1}t_n=0}}G,
$$
$n\ge1$, with $d_n:A_{n+1}\to A_n$ having $t_n$-th component equal to
$t_n$. Once again, $d_{n-1}d_n=0$ is obvious.

To prove exactness, first note that if $\Hom_\A(X,A_\bullet)$ is exact,
then for any retract $A$ of $X$ $\Hom_\A(A,A_\bullet)$ is exact as well.
Similarly if $\Hom_\A(G_i,A_\bullet)$ are exact, so is
$\Hom_\A(\coprod_iG_i,A_\bullet)$. Thus it suffices to show that
$\Hom_\A(G,A_\bullet)$ is exact for any object $G$ from $\g$. Thus suppose
given $t_n:G\to A_n$ with $d_{n-1}t_n=0$. Then $t_n=d_nt_{n+1}$, where
$t_{n+1}:G\to A_{n+1}$ is the canonical inclusion of the $t_n$-th
component into the coproduct.
\end{proof}

\begin{Lemma}\label{1-comparison}
Suppose given two $A$-augmented chain complexes $\epsilon:A_\bullet\to A$
and $\epsilon':A'_\bullet\to A$. If $A_n$ are in $\a$ for $n\ge0$ and
$A'_\bullet$ is $\a$-exact, then there exists a chain map $f:A_\bullet\to
A'_\bullet$ over $A$ (i.~e. with $f_{-1}$ equal to the identity of $A$).
Moreover this map is unique up to a chain homotopy over $A$, i.~e. for any
two $f,f':A_\bullet\to A'_\bullet$ over $A$ there is a chain homotopy
$h_\bullet$ from $f$ to $f'$ over $A$ (which means $h_0=0$).
\end{Lemma}

\begin{proof}
Since $A_0$ is in $\a$, by $\a$-exactness of $A'_\bullet$ the map
$\Hom_\A(A_0,\epsilon')$ is surjective; in particular, there is a morphism
$f_0:A_0\to A'_0$ with $\epsilon'f_0=\epsilon$. Next, as $A_1$ is also in
$\a$, and $\epsilon'f_0d_0=\epsilon d_0=0$, again by $\a$-exactness of
$A'_\bullet$ there is a map $f_1:A_1\to A'_1$ with $f_0d_0=d'_0f_1$.
Continuing this way, one obtains a sequence of maps $f_n:A_n\to A'_n$ with
$d'_nf_n=f_{n-1}d_n$ for all $n\ge0$.

Suppose now given two such sequences $f$, $f'$. Take $h_0=0:A\to A'_0$.
Since $\epsilon'(f_0-f'_0)=0$, there is a $h_1:A_0\to A'_1$ with
$f_0-f'_0=d'_0h_1=d_0'h_1+h_0d_0$. Next since
$d'_0(f_1-f'_1-h_1d_0)=(f_0-f'_0)d_0-d'_0h_1d_0=0$, there is a $h_2:A_1\to
A'_2$ with $f_1-f'_1-h_1d_0=d'_1h_2$. Continuing one obtains the desired
chain homotopy $h$.
\end{proof}

As an immediate corollary we obtain that any two $\a$-resolutions
$A_\bullet$, $A'_\bullet$ of an object are chain homotopy equivalent,
i.~e. there are maps $f:A'_\bullet\to A_\bullet$, $f':A_\bullet\to
A'_\bullet$ with $ff'$ and $f'f$ chain homotopic to identity maps. We thus
see that all the standard ingredients for doing homological algebra are
available. So we define

\begin{Definition}
\emph{$\a$-relative left derived functors} $\L^\a_nF$, $n\ge0$, of a
functor $F:\A\to{\mathscr A}$ from $\A$ to an abelian category $\mathscr
A$ are defined by
$$
(\L^\a_nF)A=H_n(F(A_\bullet)),
$$
where $A_\bullet$ is given by any $\a$-resolution of $A$. Similarly,
$\a$-relative right derived functors of a contravariant functor
$F:\A\op\to{\mathscr A}$ are given by
$$
(\R_\a^nF)A=H^n(F(A_\bullet)).
$$
\end{Definition}

By the above lemmas, $\L^\a_nF$ and $\R_\a^nF$ are indeed functors and do
not depend on the choice of resolutions. Note also that these
constructions are functorial in $F$, i.~e. a natural transformation $F\to
F'$ induces natural transformations between the corresponding derived
functors.

In particular, we have $\a$-relative Ext-groups given by
$$
\Ext^n_\a(A,X)=\left(\R_\a^n\left(\Hom_\A(\_,X)\right)\right)A=H^n(\Hom_\A(A_\bullet,X)),
$$
for objects $A$, $X$ of $\A$ and an $\a$-exact $\a$-resolution $A_\bullet$
of $A$. Note that these groups can be equipped with the \emph{Yoneda
product}
$$
\Ext^m_\a(Y,Z)\ox\Ext^n_\a(X,Y)\to\Ext^{m+n}_\a(X,Z).
$$
On representing cocycles this product can be defined as follows: given
$\a$-exact $\a$-resolutions $X_\bullet$ of $X$ and $Y_\bullet$ of $Y$, we
can represent elements of the Ext groups in question by maps $f:Y_m\to Z$
with $fd_m=0$ and $g:X_n\to Y$ with $gd_n=0$. Then similarly to the proof
of \ref{1-comparison}, we can find maps $h_0:X_{n+1}\to Y_0$, ...,
$h_{m-1}:X_{n+m}\to Y_{m-1}$, $h_m:X_{n+m+1}\to Y_m$ giving a map of
complexes, and define $[f][g]=[fh_m]$. Standard homological algebra
argument then shows that this product is well-defined, bilinear and
associative.

\begin{Examples}

\noindent1. A typical situation for the above is given by a ringoid $\g$,
with $\A$ being the category of $\g$-modules, i.~e. of linear functors
from $\g$ to abelian groups. The abelian version of the Yoneda embedding
identifies $\g$ with the full subcategory of $\A$ with objects the
representable functors. The natural choice for $\a$ is then either the
category of \emph{free $\g$-modules}, which is the closure of this full
subcategory $\g\subset\A$ under arbitrary coproducts, or that of
\emph{projective $\g$-modules} --- the closure under both coproducts and
retracts. In particular, when $\g$ has only one object, we obtain the
classical setup for homological algebra given by a ring $R$, with $\A$
being the category of $R$-modules and $\a$ that of free or projective
$R$-modules.

\noindent2. When $\A$ has finite limits, we obtain the additive case of
derived functors from \cite{Tierney&Vogel}.
\end{Examples}

\begin{Remark}
There is an obvious dual version of the above which one obtains by
replacing $\A$ with the opposite category $\A\op$. Explicitly, chain
complexes get replaced by cochain complexes (with differentials having
degree +1 rather than -1); exactness of the complex $\Hom_\A(X,A_\bullet)$
becomes replaced by that of $\Hom_\A(A^\bullet,X)$, etc.
\end{Remark}

\begin{Example}\label{staho}
Let $\A$ be the stable homotopy category of spectra and let $\a\subset\A$
be the full subcategory consisting of finite products of Eilenberg-Mac
Lane spectra over a fixed prime field $\F_p$. Let $\mathscr A$ be the mod
$p$ Steenrod algebra. The mod $p$ cohomology functor restricted to $\a$
yields an equivalence of categories for which the following diagram
commutes
$$
\xymatrix{
\A\op\ar[r]^-{H^*}&{\mathscr A}{\mathrm-}\mathbf{Mod}=\A_{\mathscr A}\\
\a\op\ar@{^(->}[u]\ar[r]^-\sim
&{\mathscr A}{\mathrm-}\mathbf{mod}=\a_{\mathscr A}.\ar@{^(->}[u]
}
$$
Here $\A\op$ denotes the opposite category of $\A$, ${\mathscr
A}{\mathrm-}\mathbf{Mod}$ is the category of positively graded
$\A$-modules and ${\mathscr A}{\mathrm-}\mathbf{mod}$ is its full
subcategory of finitely generated free modules. Given a spectrum $X$, its
$\a$-coresolution $(A_\bullet^\epsilon,d)$
$$
...\ot A_1\ot A_0\ot X\ot0\ot0\ot...
$$
is an $X$-coaugmented chain complex in $\A$, with $A_n$ in $\a$ for
$n\ge0$, which is $\a$-coexact, that is $\Hom_\A(A_\bullet^\epsilon,A')$
is acyclic for all $A'\in\a$. Hence $(A_\bullet^\epsilon,d)$ is an
$\a\op$-resolution of $X$ in $\A\op$ which is carried by the cohomology
functor $H^*$ to an $\a_{\mathscr A}$-resolution of $H^*(X)$ in
$\A_{\mathscr A}$ above. For this reason we get for a spectrum $Y$ the
binatural equation
$$
\Ext^m_{\a\op}(X,Y)=\Ext^m_{\a_{\mathscr A}}(H^*(X),H^*(Y)).
$$
Here the left hand side $\Ext^m_{\a\op}(X,Y)$ is defined in the additive
category $\A\op$ which is the opposite of the stable homotopy category.
Moreover the right hand side is the classical Ext group
$$
\Ext^m_{\a_{\mathscr A}}(H^*(X),H^*(Y))=\Ext^m_{\mathscr A}(H^*(X),H^*(Y)).
$$
\end{Example}

\section{Secondary resolutions}

We have seen in section \ref{primary} that resolutions yield the notion of
derived functors. We now introduce secondary resolutions from which we
deduce secondary derived functors. For this we need the notion of tracks.

Recall that a \emph{track category} is a category enriched in groupoids;
in particular, for all of its objects $X$, $Y$ their hom-groupoid
$\hog{X,Y}$ is given, whose objects are maps $f:X\to Y$ and whose
morphisms, denoted $\alpha:f\then f'$, are called tracks.

Equivalently, a track category is a 2-category all of whose 2-cells are
invertible. For a track $\alpha:f\then f'$ above and maps $g:Y\to Y'$,
$e:X'\to X$, the resulting composite tracks will be denoted by
$g\alpha:gf\then gf'$ and $\alpha e:fe\then f'e$. Moreover there is a
vertical composition of tracks, i.~e. composition of morphisms in the
groupoids $\hog{X,Y}$; for $\alpha:f\then f'$ and $\beta:f'\then f''$, it
will be denoted $\beta\comp\alpha:f\then f''$. An inverse of a track
$\alpha$ with respect to this composition will be denoted $\alpha\binv$.

By a \emph{track functor} we will mean a groupoid enriched functor between
track categories.

A track category $\B$ will be also depicted as $\B_1\toto\B_0$. Here
$\B_0$ being the underlying ordinary category of $\B$ obtained by
forgetting about the tracks, whereas $\B_1$ is another ordinary category
with the same objects but with morphisms from $X$ to $Y$ being tracks
$\alpha:f\then f'$ with $f,f':X\to Y$ in $\B_0$, composite of $\alpha$ and
$\beta$ in the diagram
$$
\xymat{
Z&Y\ltwocell_f^{f'}{^\alpha}&X\ltwocell_g^{g'}{^\beta}
}
$$
being
\begin{equation}\label{interchange}
\alpha\beta=\alpha g'\comp f\beta=f'\beta\comp\alpha g:fg\then f'g'.
\end{equation}
There are thus two functors $\B_1\to\B_0$ which are identity on objects
and which send a morphism $\alpha:f\then f'$ to $f$, resp. $f'$.

A track category $\B$ has the \emph{homotopy category $\B\ho$} --- an
ordinary category obtained by identifying \emph{homotopic maps}, i.~e.
maps $f,f'$ for which there exists a track $f\then f'$. It is thus the
coequalizer of $\B_1\toto\B_0$ in the category of categories.

We now assume given a track category $\B$ such that its homotopy category
is an additive category like $\A$ from section \ref{primary},
$$
\B\ho=\A,
$$
and that moreover $\B$ has a strict zero object, that is, an object * such
that for every object $X$ of $\B$, $\hog{X,*}$ and $\hog{*,X}$ are trivial
groupoids with a single morphism. It then follows that in each $\hog{X,Y}$
there is a distinguished map $0_{X,Y}$ obtained by composing the unique
maps $X\to*$ and $*\to Y$. The identity track of this map will be denoted
just by $0$. Note that $0_{X,Y}$ may also admit non-identity self-tracks;
one however has
\begin{equation}\label{zero}
0_{Y,Z}\beta=0=\alpha0_{X,Y}
\end{equation}
for any $\alpha:f\then f'$, $f,f':Y\to Z$, $\beta:g\then g'$, $g,g':X\to
Y$.

In section \ref{atrack} we introduce the notion of an ``additive track
category'' which is the most appropriate framework for secondary derived
functors and which has the properties of the track category $\B$.

\begin{Example}\label{algexample}
The most easily described example is the track category $\calli{Ch}_\A$
whose objects are chain complexes in an additive category $\A$, maps are
chain maps, and tracks are chain homotopies.

Our basic example is the track category $\calli{Pair}_\A$; it is the full
track subcategory of $\calli{Ch}$ whose objects are chain complexes
concentrated in degrees 0 and 1 only. Thus objects $A$ of
$\calli{Pair}_\A$ are given by morphisms $\partial_A:A_1\to A_0$ in $\A$,
a map $f$ from $A$ to $B$ is a pair of morphisms $(f_1:A_1\to
B_1,f_0:A_0\to B_0)$ in $\A$ making the obvious square commute, and a
track $f\then f'$ for $f,f':A\to B$ is a morphism $\phi:A_0\to B_1$ in
$\A$ satisfying $\phi\partial_A=f_1-f'_1$ and $\partial_B\phi=f_0-f'_0$.
\end{Example}

\begin{Remark}
The secondary homology $\mathscr H$, as defined in \cite{Baues}, yields a
track functor
$$
{\mathscr H}:\calli{Ch}_{\A}\to\calli{Pair}_{\A^\Z}.
$$
Here $\A$ is an abelian category, $\A^\Z$ denotes the category of
$\Z$-graded objects in $\A$, and for a chain complex $(A,d)$ in $\A$ the
$n$-th component of ${\mathscr H}(A,d)$ is given by
$$
{\mathscr
H}_n(A,d)=\left(d_n:\textrm{Coker}(d_{n+1})\to\textrm{Ker}(d_{n-1})\right).
$$
\end{Remark}

\begin{Example}\label{topexample}
A further basic example we have in mind is the track category $\B$ which
is opposite to the category of spectra, stable maps, and tracks which are
stable homotopy classes of stable homotopies.
\end{Example}

Next we describe the secondary analogues of the notions of chain complex
and resolution in \ref{1-complex}, \ref{1-resolution}.

\begin{Definition}\label{2-complex}
A \emph{secondary chain complex} $(A,d,\d)$ in a track category $\B$ is a diagram of the
form
$$
\xymat{
\cdots\ar[r]\rrlowertwocell<-12>_{0}{}
&A_{n+2}\ar[r]|-{d_{n+1}}\rruppertwocell<12>^{0}{^\d_n}
&A_{n+1}\ar[r]|-{d_n}\rrlowertwocell<-12>_{0}{_{}{\hskip1.2em\d_{n-1}}}
&A_n\ar[r]|-{d_{n-1}}\rruppertwocell<12>^{0}{^{}}
&A_{n-1}\ar[r]
&\cdots
}
$$
i.~e. a sequence of objects $A_n$, maps $d_n:A_{n+1}\to A_n$ and tracks
$\d_n:d_nd_{n+1}\then0$, $n\in\Z$, such that for each $n$ the tracks
$$
\xymat{d_{n-1}d_nd_{n+1}\ar@{=>}[r]^{d_{n-1}\d_n}&d_{n-1}0=0}
$$
and
$$
\xymat{d_{n-1}d_nd_{n+1}\ar@{=>}[r]^{\d_{n-1}d_{n+1}}&0d_{n+1}=0}
$$
coincide. Equivalently, the track $\d_{n-1}d_{n+1}\comp d_{n-1}\d_n\binv$
in hom$_{\hog{A_{n+2},A_{n-1}}}(0,0)$ must be the identity.
\end{Definition}

It is clear that a track functor $F:\B\to\B'$ between track categories as
above (which preserves the zero object) carries a secondary chain complex
in $\B$ to a secondary chain complex in $\B'$.

\begin{Examples}\label{algcomplex}

\ 

\noindent1. In the example $\calli{Pair}_\A$, a secondary chain complex looks like
$$
\xymatrix@!{
...\ar[r]
&A_{1,n+2}\ar[r]^{d_{1,n+1}}\ar[d]|\hole^>(.7){\partial_{n+2}}
&A_{1,n+1}\ar[r]^{d_{1,n}}\ar[d]|\hole^>(.7){\partial_{n+1}}
&A_{1,n}\ar[r]^{d_{1,n-1}}\ar[d]|\hole^>(.7){\partial_n}
&A_{1,n-1}\ar[r]\ar[d]|\hole^>(.7){\partial_{n-1}}
&...\\
...\ar[r]\ar[urr]^<(.3){\d_{n+1}}
&A_{0,n+2}\ar[r]_{d_{0,n+1}}\ar[urr]^<(.3){\d_n}
&A_{0,n+1}\ar[r]_{d_{0,n}}\ar[urr]^<(.3){\d_{n-1}}
&A_{0,n}\ar[r]_{d_{0,n-1}}\ar[urr]
&A_{0,n-1}\ar[r]
&...\\
}
$$
with the equations $\partial_nd_{1,n}=d_{0,n}\partial_{n+1}$,
$d_{1,n-1}d_{1,n}=\d_{n-1}\partial_{n+1}$,
$d_{0,n-1}d_{0,n}=\partial_{n-1}\d_{n-1}$ and
$d_{1,n-1}\d_n=\d_{n-1}d_{0,n+1}$ satisfied for all $n$.

More generally for $\calli{Ch}_\A$ what one obtains is a bigraded group
$A_{m,n}$ with differentials $\partial_{m,n}:A_{m+1,n}\to A_{m,n}$,
$\partial_{m,n}\partial_{m+1,n}=0$, and maps $d_{m,n}:A_{m,n+1}\to
A_{m,n}$, $\d_{m,n}:A_{m-1,n+2}\to A_{m,n}$ satisfying the similar
equalities for all $m$ and $n$.

One thus obtains a structure strongly related to what is called
multicomplex or twisted complex in the literature; cf.
\cites{Bondal&Kapranov, Meyer, OBrian&Toledo&Tong}

\noindent2. In \cite{Takeuchi&Ulbrich}, the notion of complex of categories with
abelian group structure is investigated. One can show that a slightly
strictified version of their notion coincides with that of the secondary
chain complex in an appropriate track category. On the other hand we could
relax the requirement of existence of the strict zero object to that of a
weak zero object; then the construction of \cite{Takeuchi&Ulbrich} would
be subsumed in full generality.
\end{Examples}

\begin{Definition}\label{2-map}
A \emph{secondary chain map} $(f,\phi)$ between secondary chain complexes $(A,d,\d)$
and $(A',d',\d')$  is a sequence of maps and tracks as indicated
$$
\includegraphics{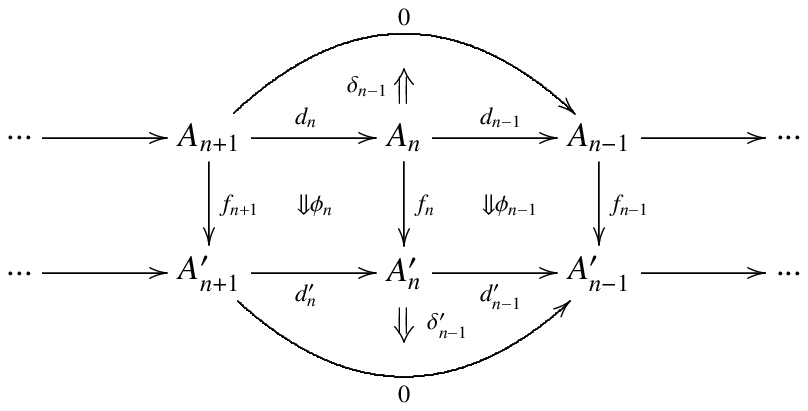}
$$
such that pasting of tracks in this diagram yields the trivial track
$0:0\then0$, that is, the resulting track diagrams
\begin{equation}\label{2-maptrack}
\xymat{
&d'_{n-1}f_nd_n\ar@{=>}[r]^{d'_{n-1}\phi_n}&d'_{n-1}d'_nf_{n+1}\ar@{=>}[dr]^{\d'_{n-1}f_{n+1}}\\
f_{n-1}d_{n-1}d_n\ar@{=>}[ur]^{\phi_{n-1}d_n}\ar@{=>}[r]^{f_{n-1}\d_{n-1}}
&f_{n-1}0\ar@{=}[r]
&0\ar@{=}[r]
&0f_{n+1}
}
\end{equation}
commute for all $n\in\Z$.

For secondary chain maps $(f,\phi):(A,d,\d)\to(A',d',\d')$ and
$(f',\phi'):(A',d',\d')\to(A'',d'',\d'')$, their \emph{composite} is given
by $(f'_nf_n,\phi'_nf_{n+1}\comp f'_n\phi_n)$, $n\in\Z$. It is
straightforward to check that this indeed defines a secondary chain map,
and that the resulting composition operation is associative. Thus these operations
determine the category of secondary chain complexes.
\end{Definition}

We now similarly to section \ref{primary} fix a full track subcategory $\b$
of $\B$, with $\a=\b\ho$.

\begin{Definition}
For a secondary complex $(A,d,\d)$ in $\B$ and an integer $n$, its
\emph{$\b$-chain of degree $n$} is a map $X\to A_n$ for some
object $X$ of $\b$. A \emph{$\b$-cycle} is a pair
$(c,\gamma)$ consisting of a $\b$-chain $c:X\to A_n$ and a
track $\gamma:d_{n-1}c\then0$ such that the track
$d_{n-2}\gamma:d_{n-2}d_{n-1}c\then d_{n-2}0=0$ is equal to
$\d_{n-2}c:d_{n-2}d_{n-1}c\then 0c=0$. A  $\b$-cycle
$(b,\beta)$ of degree $n$ is a \emph{$\b$-boundary} if there exists a $\b$-chain $a$
of degree $n+1$ and a track $\alpha:b\then d_na$ such that the following diagram of tracks
commutes:
$$
\xymat{
&d_{n-1}d_na\ar@{=>}[r]^{\d_{n-1}a}&0a\ar@{=}[dr]\\
d_{n-1}b\ar@{=>}[rrr]^\beta\ar@{=>}[ur]^{d_{n-1}\alpha}&&&0.
}
$$
A secondary complex $(A,d,\d)$ is called \emph{$\b$-exact} if all of its
$\b$-cycles are $\b$-boundaries. In other words, every diagram consisting
of solid arrows below
$$
\xymatrix@!C=5em{
&
&X\ar[dr]^0\ar[d]^c\ar@{-->}[dl]_a\ar@{*{}*{}*{}}@/_3ex/[dr]|<(.3){\buildrel\gamma\over\Longrightarrow}\\
\cdots\ar[r]
&A_{n+1}\ar@{*{}*{}*{}}@/_2.5ex/[ur]|>(.7){\buildrel\alpha\over{==\then}}
\ar[r]|{d_n}\ar@/_7ex/[rr]_0\rrlowertwocell\omit{<3>{\d_{n-1}\hskip3em}}
&A_n\ar[r]|{d_{n-1}}\ar@/_7ex/[rr]_0|>(.265)\hole\rrlowertwocell\omit{<3>{\hskip1.3em\d_{n-2}}}
&A_{n-1}\ar[r]|{d_{n-2}}
&A_{n-2}\ar[r]
&\cdots
}
$$
in which the pasted track from $d_{n-2}0_{X,A_{n-1}}$ to
$0_{A_n,A_{n-2}}c$ is the identity track
can be completed by the dashed arrows in such a way that the resulting
pasted track from $0_{A_{n+1},A_{n-1}}a$ to $0_{X,A_{n-1}}$ is the identity track.
\end{Definition}

\begin{Example}\label{total}
Consider the track category $\calli{Pair}_\A$ from \ref{algexample}, with
$\A$ the category of modules over a ring $R$, and choose
for $\b$ the full track subcategory on the objects $0\to R^n$, $n\ge0$.
Then for a secondary chain complex as in \ref{algcomplex}, a secondary
cycle of degree $n$ is a pair $(c,\gamma)\in A_{0,n}\x A_{1,n-1}$
satisfying $d_{0,n-1}c=\partial_{n-1}\gamma$ and
$\d_{n-2}c=d_{1,n-2}\gamma$. Such a cycle is a boundary if there exist
elements $a\in A_{0,n+1}$ and $\alpha\in A_{1,n}$ with
$c=d_{0,n}a+\partial_n\alpha$ and $\gamma=\d_{n-1}a+d_{1,n-1}\alpha$.

Note that we can arrange for a \emph{total complex}
$$
...\ot
A_{0,n-1}\oplus A_{1,n-2}
\xot{
\left(
\begin{smallmatrix}
d_{0,n-1}&-\partial_{n-1}\\
\d_{n-2}&-d_{1,n-2}
\end{smallmatrix}
\right)
}
A_{0,n}\oplus A_{1,n-1}
\xot{
\left(
\begin{smallmatrix}
d_{0,n}&-\partial_n\\
\d_{n-1}&-d_{1,n-1}
\end{smallmatrix}
\right)
}
A_{0,n+1}\oplus A_{1,n}
\ot
...
$$
in such a way that secondary cycles and boundaries will become usual cycles
and boundaries in this total complex. In particular then, secondary
exactness of the secondary chain complex of type \ref{algcomplex} is
equivalent to the exactness in the ordinary sense of the above total
complex.
\end{Example}

We now turn to the secondary analog of the notion of resolution from
\ref{1-resolution}.

\begin{Definition}
For an object $B$ in $\B$, a \emph{$B$-augmented secondary chain complex}
is a secondary chain complex $(B,d,\d)$ with $B_{-1}=B$, $B_{-n}=0$ for
$n>1$, and $\d_{-n}$ equal to identity track for $n>1$. For a full track
subcategory $\b$ of $\B$, a $B$-augmented secondary chain complex is
called a \emph{$\b$-resolution} of $B$ if it is $\b$-exact as a secondary
chain complex and moreover all $B_n$ for $n\ge0$ belong to $\b$.

As in the primary case, denoting $\epsilon=d_{-1}$,
$\hat\epsilon=\d_{-1}$, a $B$-augmented secondary chain complex can be considered as a
secondary chain map $(\epsilon,\hat\epsilon):B_\bullet\to B$ from the
secondary chain complex $B_\bullet$ given by $...\to B_1\to
B_0\to0\to0\to...$ with $\d_{-n}$ identities for all $n>0$, to the
secondary chain complex $B$ concentrated in degree 0, with trivial
differentials:
$$
\xymat{
\dots\ar[r]
&B_2\ar[r]^{d_1}\rruppertwocell<12>^{0}{^\hskip-2ex\d_0}\ar[d]
&B_1\ar[r]^{d_0}\ar[d]\ar@{}[dl]|{=}
&B_0\ar[r]\ar[d]^{\epsilon}\ar@{}[dl]|{\ \
\Downarrow{\scriptstyle{\hat\epsilon}}}
&0\ar[r]\ar[d]\ar@{}[dl]|{=}
&\dots\\
\dots\ar[r]
&0\ar[r]
&0\ar[r]
&B\ar[r]
&0\ar[r]
&\dots.\\
}
$$
Accordingly such an augmented secondary chain complex will be denoted
$B_\bullet^{\hat\epsilon,\epsilon}$ , and the pair
$(\epsilon,\hat\epsilon)$ will be called its augmentation.

Dually, we have the notion of a \emph{$B$-coaugmented secondary chain
complex} - the one satisfying $B_1=B$, $B_n=0$ for $n>1$, and $\d_n$ equal
to the identity track for $n>1$. Accordingly, there is a notion of a
$\b$-coresolution of $B$.
\end{Definition}

To have the analog of \ref{1-existence} we need appropriate notion of
coproduct; we might in principle use groupoid enriched, or \emph{strong}
coproducts, but for further applications more suitable is the less
restrictive notion of \emph{weak coproduct}, which we now recall.

\begin{Definition}
A family of maps $(i_k:A_k\to A)_{k\in K}$ in a track category is a
\emph{weak} (respectively, \emph{strong}) \emph{coproduct diagram} if for
every object $X$ the induced functor 
$$
\hog{A,X}\to\prod_{k\in K}\hog{A_k,X}
$$ 
is an equivalence (resp., isomorphism) of groupoids.

Thus being a weak coproduct diagram means two things:
\begin{itemize}
\item[1)] for any object $X$ and any maps $x_k:A_k\to X$, $k\in K$, there is a map $x:A\to X$
and a family of tracks $\chi_k:x_k\then xi_k$, $k\in K$;
\item[2)] for any $x,x':A\to X$ and any family of tracks $(\chi_k:xi_k\then
x'i_k)_{k\in K}$ there is a unique track $\chi:x\to x'$ satisfying $\chi_k=\chi i_k$ for all
$k\in K$,
\end{itemize}
whereas for a strong coproduct one must have
\begin{itemize}
\item[1')] for any object $X$ and any maps $x_k:A_k\to X$, $k\in K$, there
is a unique map $x:A\to X$ satisfying $x_k=xi_k$ for all $k\in K$
\end{itemize}
and 2).
\end{Definition} 

We can also weaken the notion of retract in \ref{1-existence}: call an
object $X$ a \emph{weak retract} of an object $Y$ if there exist maps
$j:X\to Y$, $p:Y\to X$ and a track $1_X\then pj$.

\begin{Lemma}\label{2-existence}
Suppose
\begin{itemize}
\item
weak coproducts of any families of objects of $\b$ exist in $\B$ and belong
to $\b$ again;
\item
there is a small track subcategory $\g$ of $\b$ such that every object of $\b$ is a
weak retract of a weak coproduct of a family of objects from $\g$.
\end{itemize}
Then every object of $\B$ has a $\b$-resolution.
\end{Lemma}

\begin{proof}
The first step is exactly as in the primary case: for an object $B$ we take
$$
B_0=\coprod_{\substack{G\in\g\\b:G\to B}}G,
$$
where $\coprod$ denotes weak coproduct. Thus in particular there is a map
$d_{-1}:B_0\to B$ and a family of tracks $\iota_b:b\then d_{-1}i_b$ for
each $b:G\to B$.

Suppose now given a $-1$-dimensional $\b$-cycle $(b,\beta)$ in the
resolution. This means just a map $b:X\to B$ for an object $X$ of $\b$,
since $\beta:d_{-2}b\then0$ is then necessarily the trivial track. By
hypothesis we then can find some weak coproduct $G=\coprod_{k\in K}G_k$ of
objects from $\g$, maps $j:X\to G$ and $p:G\to X$, and a track
$\theta:1_X\then pj$. Then by the weak coproduct property, for the maps
$i_{bpi_k}:G_k\to B_0$, where $i_k:G_k\to G$ are the weak coproduct
inclusions, there exists a map $f_0:G\to B_0$ and a family of tracks
$\iota_k:i_{bpi_k}\then f_0i_k$, $k\in K$. This then gives composite
tracks
$$
\xymatrix@1{
d_{-1}f_0i_k\ar@{<=}[r]^{d_{-1}\iota_k}
&d_{-1}i_{bpi_k}\ar@{<=}[r]^{\iota_{bpi_k}}
&bpi_k.
}
$$
Again by the weak coproduct universality there is then a track
$\phi:bp\then d_{-1}f_0$ with $d_{-1}\iota_k\comp\iota_{bpi_k}=\phi i_k$ for
all $k\in K$. Denoting $f_0j$ by $a$, one then obtains a track
$\alpha:b\then d_{-1}a$, namely the composite
$$
\xymatrix@1{b\ar@{=>}[r]^{b\theta}&bpj\ar@{=>}[r]^{\phi j}&d_{-1}f_0j,}
$$
which means that $(b,\beta)$ is a boundary, since both $\beta$ and
$\d_{-2}a\comp d_{-2}\alpha$ are zero by trivial reasons.

We next take
$$
B_1=\coprod_{\substack{G\in\g\\t_0:G\to B_0\\\tau:d_{-1}t_0\then0}}G.
$$
Then by the weak coproduct property, for the family $(t_0:G\to
B_0)_{\tau:d_{-1}t_0\then0}$ there exists $d_0:B_1\to B_0$ and tracks
$\iota_\tau:t_0\then d_0i_\tau$, where the $i_\tau:G\to B_1$ are the
coproduct structure maps. Moreover for the family
$$
\left(
\xymatrix@1{
d_{-1}d_0i_\tau\ar@{=>}[r]^{d_{-1}\iota_\tau\binv}
&d_{-1}t_0\ar@{=>}[r]^-\tau
&0=0_{G,B}=0_{B_1,B}i_\tau
}
\right)_{\tau:d_{-1}t_0\then0}
$$
there exists $\d_{-1}:d_{-1}d_0\then0$ with
$$
\d_{-1}i_\tau=\tau\comp d_{-1}\iota_\tau\binv
\leqno{(\dag_0)}
$$
for all $\tau:d_{-1}t_0\then0$. Since $\d_{-2}$ by definition must be the
identity track of the zero map, whereas $d_{-2}$ is the unique map to the
zero object, the condition $d_{n-1}\d_n=\d_{n-1}d_{n+1}$ from
\ref{2-complex} is trivially satisfied at $n=-1$.

To prove exactness at $B_0$, suppose given $b_0:X\to B_0$ and
$\beta:d_{-1}b_0\then0$, for some object $X$ of $\b$. By hypothesis, there
is a weak retraction $j:X\to G$, $p:G\to X$, $\theta:1_X\then pj$ for some
weak coproduct $G=\coprod_{k\in K}G_k$ of objects $G_k$ from $\g$. Then
for the family $(i_{\beta pi_k}:G_k\to B_1)_{k\in K}$, where $i_k:G_k\to
G$ are the coproduct structure maps, there exists a map $f_1:G\to B_1$ and
tracks $\iota_k:i_{\beta pi_k}\then f_1i_k$, $k\in K$. One thus obtains
the composite tracks
$$
\xymatrix@1{
d_0f_1i_k\ar@{<=}[r]^{d_0\iota_k}
&d_0i_{\beta pi_k}\ar@{<=}[r]^-{\iota_{\beta pi_k}}
&b_0pi_k.
}
$$
Then again by the second property of weak coproducts there is a track
$\phi_0:b_0p\then d_0f_1$ with $\phi_0i_k=d_0\iota_k\comp\iota_{\beta
pi_k}$, $k\in K$. One then gets $a_1=f_1j$ and $\alpha=\phi_0j\comp
b_0\theta:b_0\then d_0a_1$. To prove that $(a_1,\alpha)$ exhibits
$(b_0,\beta)$ as a boundary, it remains to show
$\beta=\d_{-1}a_1\comp d_{-1}\alpha$, i.~e.
$\beta=\d_{-1}f_1j\comp d_{-1}\phi_0j\comp d_{-1}b_0\theta$. Now we
have
$$
\alignbox{
\d_{-1}f_1i_k\comp d_{-1}\phi_0i_k
&=\d_{-1}f_1i_k\comp d_{-1}d_0\iota_k\comp d_{-1}\iota_{\beta pi_k}\\
&=\d_{-1}i_{\beta pi_k}\comp d_{-1}\iota_{\beta pi_k}.
}
$$
On the other hand by ($\dag_0$) one has $\d_{-1}i_{\beta pi_k}=\beta pi_k\comp
d_{-1}\iota_{\beta pi_k}\binv$, so one obtains
$$
\d_{-1}f_1i_k\comp d_{-1}\phi_0i_k=\beta pi_k
$$
for all $k$; by the weak coproduct property this then implies
$\d_{-1}f_1\comp d_{-1}\phi_0=\beta p$, hence 
$$
\alignbox{
\d_{-1}f_1j\comp d_{-1}\phi_0j\comp d_{-1}b_0\theta
&=\beta pj\comp d_{-1}b_0\theta\\
&=0\theta\comp\beta\\
&=\beta.
}
$$

Now take some $n\ge1$ and suppose all the $B_i$, $d_{i-1}$ and $\d_{i-2}$
have been already constructed for $i\le n$ in such a way that the
conditions of \ref{2-complex} and $\b$-exactness are satisfied up to
dimension $n-1$. Moreover we can assume by induction that exactness is
\emph{constructively established for $\b$-cycles originating at $\g$},
that is, for each $(n-1)$-cycle $(t_{n-1}:G\to
B_{n-1},\tau_{n-2}:d_{n-2}t_{n-1}\then0)$, $G\in\g$, with
$d_{n-3}\tau_{n-2}=\d_{n-3}t_{n-1}$, we are given explicit maps
$i_{\tau_{n-2}}:G\to B_{n-1}$ and tracks $\iota_{\tau_{n-2}}:t_{n-1}\then
d_{n-1}i_{\tau_{n-2}}$ satisfying $\tau_{n-2}=\d_{n-2}i_{\tau_{n-2}}\comp
d_{n-2}\iota_{\tau_{n-2}}$. At least this induction hypothesis is
certainly satisfied for $n=1$, by ($\dag_0$) above.

We then define
$$
B_{n+1}=\coprod_{\substack{G\in\g\\t_n:G\to
B_n\\\tau_{n-1}:d_{n-1}t_n\then0\\d_{n-2}\tau_{n-1}=\d_{n-2}t_n}}G.
$$
Then for the family $(t_n:G\to
B_n)_{\setof{\tau_{n-1}:d_{n-1}t_n\then0}{d_{n-2}\tau_{n-1}=\d_{n-2}t_n}}$
there exists $d_n:B_{n+1}\to B_n$ and tracks $\iota_{\tau_{n-1}}:t_n\then
d_ni_{\tau_{n-1}}$, where the $i_{\tau_{n-1}}:G\to B_{n+1}$ are the
coproduct structure maps. Moreover for the family
$$
\left(
\xymatrix@1{
d_{n-1}d_ni_{\tau_{n-1}}\ar@{=>}[r]^{d_{n-1}\iota_{\tau_{n-1}}\binv}
&d_{n-1}t_n\ar@{=>}[r]^-{\tau_{n-1}}
&0=0i_{\tau_{n-1}}
}
\right)_{\setof{\tau_{n-1}:d_{n-1}t_n\then0}{d_{n-2}\tau_{n-1}=\d_{n-2}t_n}}
$$
there exists $\d_{n-1}:d_{n-1}d_n\then0$ with
$$
\d_{n-1}i_{\tau_{n-1}}=\tau_{n-1}\comp d_{n-1}\iota_{\tau_{n-1}}\binv
\leqno{(\dag_n)}
$$
for all $\tau_{n-1}:d_{n-1}t_n\then0$ with $d_{n-2}\tau_{n-1}=\d_{n-2}t_n$.
To prove the condition $d_{n-2}\d_{n-1}=\d_{n-2}d_n$ from \ref{2-complex},
it suffices by the weak coproduct property to prove
$$
d_{n-2}\tau_{n-1}i_{\tau_{n-1}}\comp
d_{n-2}d_{n-1}\iota_{\tau_{n-1}}=\d_{n-2}d_ni_{\tau_{n-1}}\comp
d_{n-2}d_{n-1}\iota_{\tau_{n-1}}
$$
for each $\tau_{n-1}:d_{n-1}t_n\then0$, $t_n:G\to B_n$, with
$d_{n-2}\tau_{n-1}=\d_{n-2}t_n$. Now by ($\dag_n$) we have
$$
d_{n-2}\d_{n-1}i_{\tau_{n-1}}\comp
d_{n-2}d_{n-1}\iota_{\tau_{n-1}}=d_{n-2}\tau_{n-1},
$$
whereas by naturality we have 
$$
\d_{n-2}d_ni_{\tau_{n-1}}\comp
d_{n-2}d_{n-1}\iota_{\tau_{n-1}}=0\iota_{\tau_{n-1}}\comp\d_{n-2}t_n=\d_{n-2}t_n.
$$

Next note that the maps $i_{\tau_{n-1}}$ and tracks $\iota_{\tau_{n-1}}$
fulfil the induction hypothesis, i.~e. explicitly exhibit cycles with
domains from $\g$ as boundaries. Finally to prove exactness at $B_n$,
suppose given any $X$, any weak coproduct $G=\coprod_{k\in K}G_k$ of
objects from $\g$, any weak retraction $j:X\to G$, $p:G\to X$,
$\theta:1_X\then pj$, and any $b_n:X\to B_n$,
$\beta_{n-1}:d_{n-1}b_n\then0$ with $d_{n-2}\beta_{n-1}=\d_{n-2}b_n$. Then
for each coproduct inclusion $i_k:G_k\to G$ one has cycles given by
$b_npi_k:G_k\to B_n$, $\beta_{n-1}pi_k:d_{n-1}b_npi_k\then0$, hence for
the family $(i_{\beta_{n-1}pi_k}:G_k\to B_{n+1})_{k\in K}$ there exists a
map $f_{n+1}:G\to B_{n+1}$ and tracks $\iota_k:i_{\beta_{n-1}pi_k}\then
f_{n+1}i_k$. We then can consider the composite tracks
$$
\xymatrix@1{
d_nf_{n+1}i_k\ar@{<=}[r]^{d_n\iota_k}
&d_ni_{\beta_{n-1}pi_k}\ar@{<=}[r]^-{\iota_{\beta_{n-1}pi_k}}
&b_npi_k
}
$$
and by the weak coproduct property of $G$ find for them a track
$\phi_n:b_np\then d_nf_{n+1}$ with
$d_n\iota_k\comp\iota_{\beta_{n-1}pi_k}=\phi_ni_k$ for all
$k\in K$. This gives us an $(n+1)$-chain $a_{n+1}=f_{n+1}j$ and a track
$\alpha=\phi_nj\comp b_n\theta:b_n\then d_na_{n+1}$. To show that these
exhibit $(b_n,\beta_{n-1})$ as a boundary, one has to prove
$\beta_{n-1}=\d_{n-1}a_{n+1}\comp d_{n-1}\alpha$. The proof goes exactly as
for the case $n=0$ above.
\end{proof}

Now to the analog of \ref{1-comparison}.

\begin{Lemma}\label{2-comparison}
Suppose given two $B$-augmented secondary chain complexes $B_\bullet$ and $B'_\bullet$.
If all $B_n$ belong to $\b$ and $B'_\bullet$ is $\b$-exact, then there exists a secondary chain map
$(f,\phi):B_\bullet\to B'_\bullet$ over $B$ (i.~e. with $f_{-1}$ equal to
the identity of $B$).
\end{Lemma}

\begin{proof}
The pair $d_{-1}:B_0\to B$, identity$_0:d_{-2}d_{-1}\then0$  can be
considered as a $(-1)$-cycle in $B'_\bullet$, so by $\b$-exactness of
$B'_\bullet$ there exist $f_0:B_0\to B'_0$ and $\phi_{-1}:d_{-1}\then
d'_{-1}f_0$. Next $f_0d_0$, $\d_{-1}\comp\phi_{-1}\binv
d_0:d'_{-1}f_0d_0\then d_{-1}d_0\then0$ is a 0-cycle in $B'_\bullet$, so
again by exactness of $B'_\bullet$ there are $f_1:B_1\to B'_1$ and
$\phi_0:f_0d_0\then d'_0f_1$ with
$$
\d_{-1}\comp\phi_{-1}\binv d_0=\d'_{-1}f_1\comp d'_{-1}\phi_0,
\leqno{(*)}
$$
which ensures the condition of \ref{2-map} for $n=0$. Then $f_1d_1$,
$f_0\d_0\comp\phi_0\binv d_1:d'_0f_1d_1\then f_0d_0d_1\then f_00=0$ is a
1-cycle in $B'_\bullet$. Indeed $(*)$ above implies
$\d'_{-1}f_1d_1=\d_{-1}d_1\comp\phi_{-1}\binv d_0d_1\comp
d'_{-1}\phi_0\binv d_1$; on the other hand $\d_{-1}d_1\comp\phi_{-1}\binv
d_0d_1=d_{-1}\d_0\comp\phi_{-1}\binv d_0d_1=\phi_{-1}\binv0\comp
d'_{-1}f_0\d_0=d'_{-1}f_0\d_0$, so $\d'_{-1}f_1d_1=d'_{-1}f_0\d_0\comp
d'_{-1}\phi_0\binv d_1$, which precisely means that the cycle condition is
fulfilled. One thus obtains $f_2:B_2\to B'_2$ and $\phi_1:f_1d_1\then
d'_1f_2$ such that $f_0\d_0\comp\phi_0\binv d_1=\d'_0f_2\comp d'_0\phi_1$,
so the condition of \ref{2-map} at $n=1$ is also satisfied.

It is clear that continuing in this way one indeed obtains a secondary
chain map.
\end{proof}

\section{Additive track categories}\label{atrack}

The secondary analogue of an additive category is an additive track
category considered in this section. For related conditions, see
\cites{Baues&Pirashvili}.

\begin{Definition}\label{tradd}
A track category $\B$ is called \emph{additive} if it has a strict zero
object $*$, the homotopy category $\A=\B\ho$ is additive and moreover $\B$
is a linear track extension
$$
D\to\B_1\toto\B_0\to\A
$$
of $\A$ by a biadditive bifunctor
$$
D:\A\op\x\A\to\calli{Ab}.
$$
Explicitly, this means the following: a biadditive bifunctor $D$ as above
is given together with a system of isomorphisms
\begin{equation}\label{action}
\sigma_f:D(X,Y)\to\Aut_{\hog{X,Y}}(f)
\end{equation}
for each 1-arrow $f:X\to Y$ in $\B$, such that for any $f:X\to Y$, $g:Y\to
Z$, $a\in D(X,Y)$, $b\in D(Y,Z)$, $\alpha:f\then f'$ one has
$$
\begin{aligned}
\sigma_{gf}(ga)&=g\sigma_f(a);\\
\sigma_{gf}(bf)&=\sigma_g(b)f;\\
\alpha\comp\sigma_f(a)&=\sigma_{f'}(a)\comp\alpha.
\end{aligned}
$$

\begin{Remark}
Using \ref{action} we can identify the bifunctor $D$ via the natural equation
$$
D(X,Y)=\Aut(0_{X,Y}),
$$
where $0=0_{X,Y}:X\to*\to Y$ is the unique morphism factoring through the
zero object.
\end{Remark}

A \emph{strict equivalence} between additive track categories $\B$, $\B'$
is a track functor $\B\to\B'$ which induces identity on $\A$ and is
compatible with the actions \ref{action} above. Thus for fixed $\A$ and
$D$ as above, one obtains a category whose objects are additive track
categories which are linear track extensions of $\A$ by $D$ and morphisms
are strict equivalences. This category will be denoted by Trext($A$;$D$).
For an additive category $\A$ and a biadditive bifunctor $D$ on it, there
is a bijection
$$
\pi_0(\textrm{Trext}(\A;D))\approx H^3(\A;D),
$$
where $\pi_0(\mathbb{C})$ denotes the set of connected components of a
small category $\mathbb C$. Two additive track categories are called
\emph{equivalent} if they are in the same connected component of
Trext($\A$;$D$). Thus in particular (as shown in \cites{PirashviliI,Baues&Dreckmann})  each
additive track category $\B$ as above determines a class $\brk{\B}\in
H^3(\A;D)$. 
\end{Definition}

As shown in \cite{Pirashvili&Waldhausen}, when $\A$ is the category of free
finitely generated modules over a ring $R$ and $D$ is given by
$D(X,Y)=\Hom_R(X,B\ox_RY)$ for some $R$-$R$-bimodule $B$, there are
isomorphisms
$$
H^3(\A;D)\cong HML^3(R;B)\cong THH^3(HR;HB),
$$
where $HML^*$ denotes Mac Lane cohomology, $THH$ is the topological
Hochschild cohomology, and $HR$ and $HB$ are the Eilenberg-Mac Lane spectra
corresponding to $R$ and $B$. 

\begin{Definition}\label{sigmaomega}
An additive track category $\B$ is \emph{$\LL$-additive} if an additive
endofunctor $\LL:\A\to\A$ is given which left represents the bifunctor
$D$, i.~e. $\B$ is a linear track extension of $\A$ by the bifunctor
$$
D(X,Y)=\Hom_\A(\LL X,Y).
$$
Dually, $\B$ is \emph{$\RR$-additive} if an additive endofunctor
$\RR:\A\to\A$ is given such that $\B$ is a linear track extension of $\A$
by the bifunctor
$$
D(X,Y)=\Hom_\A(X,\RR Y).
$$
For objects $X$, $Y$ in a $\LL$-, resp. $\RR$-additive track category
$\B$ we will denote the group $\Hom_\A(\LL^mX,Y)$, resp.
$\Hom_\A(X,\RR^mY)$ by $[X,Y]^m$.
\end{Definition}

In examples from topology the functor $\LL$ is the suspension and the
functor $\RR$ is the loop space, compare also \cite{Baues&JibladzeII}.

$\LL$- or $\RR$-additivity of a track category enables one to relate
secondary exactness of a secondary chain complex with exactness of the
corresponding chain complex in the homotopy category.

\begin{Lemma}\label{hosec}
Let $\B$ be a track category with the additive homotopy category
$\A=\B\ho$, let $\b$ be a full track subcategory of $\B$ and denote
$\a=\b\ho$. Suppose that one of the following conditions is satisfied:
\begin{itemize}
\item[a)]
$\B$ is $\LL$-additive and $\a$ is closed under suspensions
(i.~e. for each $X\in\a$ one has $\LL X\in\a$); or
\item[b)]
$\B$ is $\RR$-additive and the functor $\RR$ is
$\a$-exact (i.~e. for an $\a$-exact complex $A_\bullet$ in $\A$, $\RR
A_\bullet$ is also $\a$-exact).
\end{itemize}
Then for any secondary chain complex $(A,d,\d)$ in $\B$, $\a$-exactness
of its image $(A,[d])$ in $\A$ implies $\b$-exactness of $(A,d,\d)$.

If moreover $(A,d,\d)$ is bounded below, then conversely its $\b$-exactness
implies $\a$-exactness of $(A,[d])$.
\end{Lemma}

\begin{proof}
Unraveling definitions, we have that for any $a_n:X\to A_n$ with $X$ in
$\b$ and for any track $\alpha_{n-1}:d_{n-1}a_n\then0$ there exists
$a_{n+1}:X\to A_{n+1}$ and a track $\alpha_n:a_n\then d_na_{n+1}$. From
this, we have then to deduce that for $(a_n,\alpha_{n-1})$ as above with
the additional property $d_{n-2}\alpha_{n-1}=\d_{n-2}a_n$ one can
actually find $(\tilde a_{n+1},\tilde\alpha_n)$ as above with the
additional property $\alpha_{n-1}=\d_{n-1}\tilde a_{n+1}\comp
d_{n-1}\tilde\alpha_n$.

Now suppose given $a_n:X\to A_n$, $\alpha_{n-1}:d_{n-1}a_n\then0$ with
$d_{n-2}\alpha_{n-1}=\d_{n-2}a_n$, and any $a_{n+1}:X\to A_{n+1}$,
$\alpha_n:a_n\then d_na_{n+1}$. Consider then the element
$\omega_{n-1}\in\Aut(d_{n-1}a_n)$ given by the composite
$$
\xymatrix@1{
d_{n-1}a_n\ar@{=>}[r]^{d_{n-1}\alpha_n}
&d_{n-1}d_na_{n+1}\ar@{=>}[r]^{\d_{n-1}a_{n+1}}
&0a_{n+1}=0\ar@{=>}[r]^{\alpha_{n-1}\binv}
&d_{n-1}a_n.
}
$$
For this element one has $d_{n-2}\omega_{n-1}=0$. Indeed, this equality
is equivalent to the equality
$$
d_{n-2}\d_{n-1}a_{n+1}\comp d_{n-2}d_{n-1}\alpha_n=d_{n-2}\alpha_{n-1}
$$
of tracks $\Aut(0_{X,A_{n-2}})$. But $d_{n-2}\alpha_{n-1}=\d_{n-2}a_n$. Moreover
by naturality there is a commutative diagram
$$
\xymat
{
d_{n-2}d_{n-1}a_n\ar@{=>}[rr]^{d_{n-2}d_{n-1}\alpha_{n+1}}\ar@{=>}[d]_{\d_{n-2}a_n}&&d_{n-2}d_{n-1}d_na_{n+1}\ar@{=>}[d]^{\d_{n-2}d_na_{n+1}}\\
0a_n\ar@{=>}[rr]^{0\alpha_{n+1}}\ar@{=}[dr]&&0d_na_{n+1}\ar@{=}[dl]\\
&0
}
$$
showing that $\d_{n-2}a_n=\d_{n-2}d_na_{n+1}\comp
d_{n-2}d_{n-1}\alpha_{n+1}$. It thus follows that $d_{n-2}\omega_{n-1}=0$ iff one
has
$$
d_{n-2}\d_{n-1}a_{n+1}\comp
d_{n-2}d_{n-1}\alpha_{n+1}=\d_{n-2}d_na_{n+1}\comp
d_{n-2}d_{n-1}\alpha_{n+1},
$$
which is clear since $(A,d,\d)$ is a secondary complex.

Now if a) is satisfied, then there is a commutative diagram
$$
\xymat{
[\LL X,A_n]\ar[r]^{[d_{n-1}]\_}\ar[d]_\cong&[\LL
X,A_{n-1}]\ar[r]^{[d_{n-2}]\_}\ar[d]_\cong&[\LL X,A_{n-2}]\ar[d]_\cong\\
\Aut(a_n)\ar[r]^{d_{n-1}\_}&\Aut(d_{n-1}a_n)\ar[r]^{d_{n-2}\_}&\Aut(d_{n-2}d_{n-1}a_n)
}
$$
Similarly if b) holds, then one has the diagram
$$
\xymat{
[X,\RR A_n]\ar[r]^{\RR[d_{n-1}]\_}\ar[d]_\cong&[X,\RR
A_{n-1}]\ar[r]^{\RR[d_{n-2}]\_}\ar[d]_\cong&[X,\RR A_{n-2}]\ar[d]_\cong\\
\Aut(a_n)\ar[r]^{d_{n-1}\_}&\Aut(d_{n-1}a_n)\ar[r]^{d_{n-2}\_}&\Aut(d_{n-2}d_{n-1}a_n).
}
$$
In both cases, it follows that there exists $\omega_n\in\Aut(a_n)$ such that
$\omega_{n-1}=d_{n-1}\omega_n$. 

Let us then choose
\begin{align*}
\tilde a_{n+1}&=a_{n+1},\\
\tilde\alpha_n&=\alpha_n\comp\omega_n\binv.
\end{align*}
Then $\alpha_{n-1}\binv\comp\d_{n-1}a_{n+1}\comp
d_{n-1}\tilde\alpha_n=\omega_{n-1}\comp d_{n-1}\omega_n\binv$ is the
identity track of $d_{n-1}a_n$, i.~e.
$$
\alpha_{n-1}=\d_{n-1}a_{n+1}\comp d_{n-1}\tilde\alpha_n,
$$
as desired.

For the converse, by boundedness we can assume by induction that $(A,[d])$
is exact in all degrees $<n$. Let us then consider any $\a$-cycle
$[c]\in[X,A_n]$ in $(A,[d])$, choose a representative map $c:X\to A_n$ and
a track $\gamma:0\then d_{n-1}c$ and consider the composite track
$\omega=\d_{n-2}c\comp d_{n-2}\gamma$ in $\Aut(0_{X,A_{n-2}})$.

The track $d_{n-3}\omega$ is the identity track $0$ of $0_{X,A_{n-3}}$.
Indeed $d_{n-3}\d_{n-2}=\d_{n-3}d_{n-1}$ by definition of a secondary
chain complex, so $d_{n-3}\omega=\d_{n-3}d_{n-1}c\comp
d_{n-3}d_{n-2}\gamma$. Then by \eqref{interchange} for
$\d_{n-3}:d_{n-3}d_{n-2}\then0_{A_{n-1},A_{n-3}}$ and
$\gamma:0_{X,A_{n-1}}\then d_{n-1}c$ one has $\d_{n-3}d_{n-1}c\comp
d_{n-3}d_{n-2}\gamma=0_{A_{n-1},A_{n-3}}\gamma\comp\d_{n-3}0_{X,A_{n-1}}$
and by \eqref{zero} both of the constituents in the last composition are
identity tracks.

Now by induction hypothesis $(A,[d])$ is $\a$-exact in degree $n-2$, hence
if a), resp. b) holds, then the diagram
$$
\xymat{
[\LL X,A_{n-1}]\ar[r]^{[d_{n-2}]\_}\ar[d]_\cong&[\LL
X,A_{n-2}]\ar[r]^{[d_{n-3}]\_}\ar[d]_\cong&[\LL X,A_{n-3}]\ar[d]_\cong\\
\Aut(0_{X,A_{n-1}})\ar[r]^{d_{n-2}\_}&\Aut(0_{X,A_{n-2}})\ar[r]^{d_{n-3}\_}&\Aut(0_{X,A_{n-3}}),
}
$$
resp.
$$
\xymat{
[X,\RR A_{n-1}]\ar[r]^{\RR[d_{n-2}]\_}\ar[d]_\cong&[X,\RR
A_{n-2}]\ar[r]^{\RR[d_{n-3}]\_}\ar[d]_\cong&[X,\RR A_{n-3}]\ar[d]_\cong\\
\Aut(0_{X,A_{n-1}})\ar[r]^{d_{n-2}\_}&\Aut(0_{X,A_{n-2}})\ar[r]^{d_{n-3}\_}&\Aut(0_{X,A_{n-3}}).
}
$$
shows that there exists $\alpha\in\Aut(0_{X,A_{n-1}})$ such that
$\omega=d_{n-2}\alpha$. Then for $\tilde\gamma=\gamma\comp\alpha\binv$ one
has $\d_{n-2}c\comp d_{n-2}\tilde\gamma=\d_{n-2}c\comp d_{n-2}\gamma\comp
d_{n-2}\alpha=\omega\comp\omega\binv=0$, so that $(c,\tilde\gamma)$ is a
secondary cycle. Then by secondary $\b$-exactness of $(A,d,\d)$ there is a
$b:X\to A_{n+1}$ and $\beta:c\then d_nb$, so $[c]$ is the boundary of $[b]$
in $[X,(A,[d])]$. Thus $(A,[d])$ is exact in degree $n$ and we are done.
\end{proof}

\section{The secondary $\Ext$}\label{sext}

In this section we deduce from a secondary resolution a differential
defined on ``primary'' derived functors as studied in section
\ref{primary}. This differential is the analogue of the $d_2$-differential
in a spectral sequence. We use the secondary differential to define
certain ``secondary'' derived functors.

Let $\B$ be an additive track category with the additive homotopy category
$\A=\B\ho$. Let us furthermore fix a full additive subcategory $\a$ in
$\A$; it determines the full track subcategory $\b$ of $\B$ on the same objects. It
is clear that if $\b$ satisfies the conditions of \ref{2-existence}, then
$\a$ will satisfy those of \ref{1-existence}. We can then consider the
$\a$-derived functors in $\A$. In particular, the $\Ext$ groups
$\Ext^n_\a(X,Y)$ are defined for any objects $X$, $Y$ in $\B$. Moreover if
$\B$ is $\LL$-, resp. $\RR$-additive, then derived functors of the
functor $D(X,Y)=\Aut(0_{X,Y})$ are given by
$$
D^n_\a(X,Y)\cong\Ext^n_\a(\LL X,Y),
$$
resp.
$$
D^n_\a(X,Y)\cong\Ext^n_\a(X,\RR Y).
$$
We will use these isomorphisms to introduce the graded $\Ext$ groups
$\Ext^n_\a(X,Y)^m=\Ext^n_\a(\LL^mX,Y)$, resp.
$\Ext^n_\a(X,Y)^m=\Ext^n_\a(X,\RR^mY)$. Evidently if $\B$ is both $\LL$- and
$\RR$-additive, these groups coincide.

We will from now on assume in what follows that for the pair $(\B,\a)$ one
of the conditions in \ref{hosec} is satisfied, i.~e. either $\B$ is
$\LL$-additive and $\a$ is closed under $\LL$ or $\B$ is
o$\RR$-additive and $\RR$ preserves $\a$-exactness of chain
complexes in $\B\ho$; moreover in the latter case we also assume that $\a$
is closed under $\RR$.

We are going to define the \emph{secondary differential}
$$
d_{(2)}=d_{(2)}^{n,m}:\Ext^n_\a(X,Y)^m\to\Ext_\a^{n+2}(X,Y)^{m+1}.
$$
Replacing, if needed, $X$ by $\LL^mX$ (resp. $Y$ by $\RR^mY$) we
might clearly assume $m=0$ here. Moreover by \ref{2-existence} we may
suppose that a $\b$-exact $\b$-resolution
$(X_\bullet,d_\bullet,\d_\bullet)$ of $X$ is given. Then by \ref{hosec} it
determines an $\a$-exact $\a$-resolution $(X_\bullet,[d_\bullet])$ of $X$
in $\A$. Hence an element of $\Ext^n_\a(X,Y)$ gets represented by an
$n$-dimensional cocycle in that resolution, i.~e. by a homotopy class
$[c]:X_n\to Y$ with $[c][d_n]=0$. Thus we may choose a map $c\in[c]$ and a
track $\gamma:0\then cd_n$ in $\B$, as in the diagram below: 
\begin{equation}\label{d2constr}
\xymat{
X_{n+3}\ar[r]^{d_{n+2}}\rrlowertwocell<-12>_{0}{\ \ \ \ \d_{n+1}}
&X_{n+2}\ar[r]|-{d_{n+1}}\rruppertwocell<12>^{0}{^\d_n}
&X_{n+1}\ar[r]|-{d_n}\ddrtwocell<0>_0{^<-4>\gamma}
&X_n\ar[dd]^c\\
\\
&&&Y
}
\end{equation}
Then the composite track $c\d_n\comp\gamma d_{n+1}\in\Aut(0_{X_{n+2},Y})$
determines an element $\Gamma=\Gamma_{c,\gamma}$ in the group
$\Aut(0_{X_{n+2},Y})$. One then has $\Gamma d_{n+2}=0$. Indeed
$$
\alignbox{
\Gamma d_{n+2}&=(c\d_n\comp\gamma d_{n+1})d_{n+2}\\
&=c\d_nd_{n+2}\comp\gamma d_{n+1}d_{n+2}\\
&=cd_n\d_{n+1}\comp\gamma d_{n+1}d_{n+2}\\
&=\gamma0\comp0\d_{n+1}\\
&=0.
}
$$
Thus $\Gamma$ determines an $(n+2)$-cocycle in
$\Aut(0_{(X_\bullet,[d_\bullet]),Y})\cong[(X_\bullet,[d_\bullet]),Y]^1$. We then have

\begin{Theorem}\label{indep}
The above construction does not depend on the choice of $c$, $\gamma$ and
the resolution, up to coboundaries in
$[(X_\bullet,[d_\bullet]),Y]^1$; hence the assignment
$[c]\mapsto[\Gamma_{c,\gamma}]$ gives a well-defined homomorphism
$$
d_{(2)}^{n,m}:\Ext^n_\a(X,Y)^m\to\Ext_\a^{n+2}(X,Y)^{m+1}.
$$
\end{Theorem}

\begin{Remark}
Of course the above homomorphism depends on the additive track category
$\B$ in which we define the secondary resolution. In fact, $d_{(2)}$
depends only on the track subcategory $\b\set{X,Y}\subset\B$ obtained by
adding to $\b$ the objects $X$ and $Y$ and all morphisms and tracks from
$\hog{Z,X}$, $\hog{Z,Y}$ for all objects $Z$ from $\b$. We shall see in
section \ref{invariance} below that additive track categories $\B$, $\B'$
with subcategories $\b$, $\b'$ such that the track categories
$\b\set{X,Y}$ and $\b'\set{X,Y}$ are track equivalent yield the same
differential $d_{(2)}$.
\end{Remark}

If the composites $d_{(2)}^{n,m}d_{(2)}^{n-2,m-1}$ are all zero (as this
is the case for examples derived from spectral sequences), we define the
\emph{secondary $\Ext$ groups}
\begin{equation}\label{defdif}
\Ext_\b^n(X,Y)^m:=\ker(d_{(2)}^{n,m})/\im(d_{(2)}^{n-2,m-1}).
\end{equation}
This then will be, in examples, the $E_3$-term of a spectral sequence. We
point out that the secondary $\Ext$-groups are well defined and do not
depend on the choice of the secondary resolution. We shall use the
secondary $\Ext$-groups for the computation of the $E_3$-term in the Adams
spectral sequence, see \cite{E3}.

\begin{proof}
We will first show that the cocycles corresponding to $(c,\gamma)$ and
$(c,\gamma')$ for $\gamma,\gamma':0\then cd_n$ are cohomologous. Indeed
the first one is $c\d_n\comp\gamma d_{n+1}$ and the second is
$$
\alignbox{
c\d_n\comp\gamma'd_{n+1}&=c\d_n\comp\gamma d_{n+1}
\comp\gamma\binv d_{n+1}\comp\gamma'd_{n+1}\\
&=c\d_n\comp\gamma d_{n+1}\comp(\gamma\binv\comp\gamma')d_{n+1},
}
$$
so these cocycles indeed differ by the coboundary of
$\gamma\binv\comp\gamma'$. Thus we obtain a map $d_{(2)}$ from the group
of $n$-cocycles of $\Hom_\A((X_\bullet,[d_\bullet]),Y)$ to
$H^{n+2}(\Aut(0_{(X_\bullet,[d_\bullet]),Y}))\cong\Ext_\a^{n+2}(X,Y)^1$.

Next let us show that the map we just constructed is actually a
homomorphism.

To see this, let us choose maps $p_1,\nabla,p_2:Y\oplus Y\to Y$ in the
homotopy classes $([1_Y],0)$, $([1_Y],[1_Y])$, $(0,[1_Y])$ $\in$ $[Y\oplus
Y,Y]$ respectively. Thus for any two maps $c_1,c_2:X\to Y$ there is a map
$c_{1,2}:X\to Y\oplus Y$ such that there exist tracks
$\pi_i:p_ic_{1,2}\then c_i$, $i=1,2$, and moreover $[c_1]+[c_2]=[\nabla
c_{1,2}]$. Now suppose $c_i$ represent cocycles, then 
$[c_{1,2}][d_n]=([c_1][d_n],[c_2][d_n])=(0,0)\in[X_{n+1},Y]\x[X_{n+1},Y]\approx[X_{n+1},Y\oplus
Y]$, so there is a track $\gamma:0\then c_{1,2}d_n$. Consequently the
cohomology class $d_{(2)}([c_1]+[c_2])=d_{(2)}([\nabla c_{1,2}])$ can be
represented by the cocycle
$$
\nabla c_{1,2}\d_n\comp\nabla\gamma d_{n+1}=\nabla(c_{1,2}\d_n\comp\gamma d_{n+1}).
$$
On the other hand $d_{(2)}([f_i])$, $i=1,2$, can be represented by
$$
\alignbox{
c_i\d_n\comp\pi_id_nd_{n+1}\comp p_i\gamma d_{n+1}
&=\pi_i0\comp p_ic_{1,2}\d_n\comp p_i\gamma d_{n+1}\\
&=p_ic_{1,2}\d_n\comp p_i\gamma d_{n+1}\\
&=p_i(c_{1,2}\d_n\comp\gamma d_{n+1})
}
$$
(see the diagram below).
$$
\xymat{
X_{n+2}\ar[r]_-{d_{n+1}}\rruppertwocell<12>^{0}{^{\d_n\ }}
&X_{n+1}\ar[r]|-{d_n}\ddrtwocell<0>^0{^<-4>\gamma}\xtwocell[3,1]{}<0>_0{=<-1.2>}
&X_n\ar[dd]|{c_{1,2}}\ddruppertwocell<0>^{c_i}{^<4>\pi_i}\\
\\
&&Y\oplus Y\ar[r]_{p_i}\ar[d]^\nabla&Y\\
&&Y
}
$$
But by assumption $\Aut(0)$ is biadditive, which in particular means that
the map
$$
(p_1\_,p_2\_):\Aut(0_{X,Y\oplus Y})\to\Aut(0_{X,Y})\x\Aut(0_{X,Y})
$$
is an isomorphism, and moreover addition in $\Aut(0_{X,Y})$ is given by
the composite of the left action $\nabla\_$ with the inverse of that
isomorphism. This obviously means
$d_{(2)}([c_1]+[c_2])=d_{(2)}([c_1])+d_{(2)}([c_2])$.

It follows that in order to show that $d_{(2)}$ factors through a homomorphism from the
group $\Ext^n_\a(X,Y)$ $=$ $H^n([(X_\bullet,[d_\bullet]),Y])$ it suffices
to show that $d_{(2)}$ vanishes on coboundaries, i.~e. on cocycles of the
form $[c]=[ad_{n-1}]$, for some map $a:X_{n-1}\to Y$. But for such a
cocycle we may choose the track $\gamma:0\then ad_{n-1}d_n$ to be
$a\d_{n-1}\binv$, and then the value of $d_{(2)}$ on it will be
represented by the cocycle $ad_{n-1}\d_n\comp a\d_{n-1}d_n=0$ --- see the
diagram.
$$
\xymat{
X_{n+2}\ar[r]_-{d_{n+1}}\rruppertwocell<8>^0{^{\d_n\ }}
&X_{n+1}\ar[r]|-{d_n}\rrlowertwocell<-8>_0{^{\hskip4em
\d_{n-1}\binv}}\ar[ddr]_0
&X_n\ar[r]^{d_{n-1}}
&X_{n-1}\ar[ddl]^a\\
&&=\\
&&Y
}
$$

Finally we must show that $d_{(2)}$ does not depend on the choice of the
secondary resolution. Indeed consider any two $\b$-exact $\b$-resolutions 
$(X_\bullet,d_\bullet,\d_\bullet)$ and
$(X'_\bullet,d'_\bullet,\d'_\bullet)$ of $X$. By \ref{2-comparison} there
is a secondary chain map $(f,\phi)$ between them over $X$. Obviously then
$[f]$ determines a chain map between $(X_\bullet,[d_\bullet])$ and
$(X'_\bullet,[d'_\bullet])$ inducing isomorphisms $f^*$ on cohomology of
the cochain complexes obtained by applying $[\_,Y]$ and $\Aut(0_{\_,Y})$.
We must then show that the diagrams
$$
\xymat{
H^n([(X_\bullet,[d_\bullet]),Y])\ar[r]^{d_{(2)}}\ar[d]_{f^*}
&H^n(\Aut(0_{(X_\bullet,[d_\bullet]),Y}))\ar[d]^{f^*}\\
H^n([(X'_\bullet,[d'_\bullet]),Y])\ar[r]^{d_{(2)}}
&H^n(\Aut(0_{(X'_\bullet,[d'_\bullet]),Y}))
}
$$
commute. This can be seen from the diagram
$$
\xymat{
X'_{n+2}\ar[r]^{d'_{n+1}}\rruppertwocell<12>^{0}{^\hskip-2ex\d'_n}\ar[d]_{f_{n+2}}
&X'_{n+1}\ar[r]^{d'_n}\ar[d]|{f_{n+1}}
\ar@{}[dl]|{\ \ \Downarrow{\scriptstyle{\phi_{n+1}}}}
&X'_n\ar[d]^{f_n}
\ar@{}[dl]|{\ \ \Downarrow{\scriptstyle{\phi_n}}}\\
X_{n+2}\ar[r]_{d_{n+1}}\rrlowertwocell<-10.5>_{0}{_\hskip1ex\d_n}
&X_{n+1}\ar[r]_{d_n}\ar[dr]|{\hole}_>(.8)0
\ar@{}[drr]|<(.4){\buildrel{\scriptstyle{\gamma}}\over\Rightarrow}
&X_n\ar[d]^c\\
&&Y;&\ 
}
$$
in more detail, one observes the track diagram
$$
\xymat{
&cf_n0\ar@{<=}[r]^{cf_n\d'_n}\ar@{=}[dd]\ar@{=}[dl]
&cf_nd'_nd'_{n+1}\ar@{<=}[r]^{c\phi_n\binv d'_{n+1}}
&cd_nf_{n+1}d'_{n+1}\ar@{<=}[r]^{\gamma f_{n+1}d'_{n+1}}\ar@{=>}[dd]_{cd_n\phi_{n+1}}
&0f_{n+1}d'_{n+1}\ar@{=>}[dd]_{0\phi_{n+1}}\ar@{=}[dr]\\
0\ar@{=}[dr]&&&&&0\ar@{=}[dl]\\
&c0f_{n+2}\ar@{<=}[rr]^{c\d_nf_{n+2}}&&cd_nd_{n+1}f_{n+2}\ar@{<=}[r]^{\gamma
d_{n+1}f_{n+2}}&0d_{n+1}f_{n+2}
}
$$
whose left part commutes by \ref{2-maptrack} and the right part by
naturality. Now the lower composition of this diagram is
$$
(c\d_n\comp\gamma d_{n+1})f_{n+2}=f^*d_{(2)}([c]),
$$
whereas the upper one is
$$
(cf_n)\d'_n\comp(c\phi_n\binv\comp\gamma f_{n+1})d'_{n+1},
$$
which represents $d_{(2)}(f^*([c]))$, since we might choose for
$\gamma':0\then f^*(c)d'_n$ the track $c\phi_n\binv\comp\gamma
f_{n+1}:0\then cf_nd'_n=f^*(c)d'_n$.
\end{proof}

\section{Invariance of the secondary differential in the equivalence class
of the track extension}\label{invariance}

In this section we will prove that the secondary differential
$$
d_{(2)}^{n,m}:\Ext^n_\a(X,Y)^m\to\Ext_\a^{n+2}(X,Y)^{m+1}
$$
constructed from an additive track category $\B$ depends only on the class
$$
\brk{\B}\in H^3(\A;D).
$$
More precisely one has
\begin{Theorem}
Suppose given additive track categories
$\B$ and $\B'$ with $\B\ho=\A=\B'\ho$. Then for any additive subcategory $\a\subset\A$, the
secondary differentials $d^{n,m}_{(2)}$ constructed from $\B$ and $\B'$
coincide provided there is a strict equivalence of track subcategories
$\b\set{X,Y}\xto\sim\b'\set{X,Y}$.
\end{Theorem}

\begin{proof}
Recall the construction of
$$
d_{(2)}^{n,m}:\Ext^n_\a(X,Y)^m\to\Ext_\a^{n+2}(X,Y)^{m+1}.
$$
Let $\b\subset\B$, $\b'\subset\B'$ be the full track subcategories in $\B$,
resp. $\B'$, on objects from $\a$. Then $F(\b)\subset\b'$.
One starts from a $\b$-exact $\b$-resolution
$(X_\bullet,d_\bullet,\d_\bullet)$ of $X$ in $\B$; according to \ref{indep},
the resulting $d_{(2)}$ does not depend on the choice of such a resolution.
Suppose now given an element in
$\Ext_\a^n(X,Y)$ represented by a $Y$-valued $n$-cocycle $[c]\in[X_n,Y]$ in the $\a$-exact
$\a$-resolution $(X_\bullet,[d_\bullet])$ of $X$ in $\A$. By our
construction, value on this element of the $d_{(2)}$ corresponding to $\B$
is obtained by choosing a representative $[c]\ni c:X_n\to Y$ and a
track $\gamma:0\then cd_n$ in $\B$, as in \ref{d2constr}. One then has
$$
d_{(2)}([c])=[c\d_n\comp\gamma d_{n+1}].
$$

But it is clear that $F(X_\bullet,d_\bullet,\d_\bullet)$ is a $\b'$-exact
$\b'$-resolution of $X$ in $\B'$. We then might choose $F(c)$ and
$F(\gamma)$ for the corresponding data in $\B'$, which would give us the
element of $\Aut_{\B'}(0_{X_{n+2},Y})$ equal to $F(c)F(\d_n)\comp F(\gamma)F(d_{n-1})=F(c\d_n\comp\gamma
d_{n_1})$. Since by assumption $F$ induces identity on $\Aut(0)$, theorem
follows.
\end{proof}

\section{Resolutions of the Adams type}

Let $\B$ be a track category with a strict zero object and homotopy
category $\B\ho=\A$.

\begin{Definition}\label{coaugdef}
For an object $X$ of $\A$, an \emph{$X$-coaugmented sequence} $\r$ is
a diagram in $\A$ of the form
$$
\xymatrix{
...
&Y_{n+1}\ar[l]
&A_n\ar[l]_-{\p_n}
&Y_n\ar[l]_-{\i_n}
&...\ar[l]
&Y_2\ar[l]
&A_1\ar[l]_-{\p_1}
&Y_1\ar[l]_-{\i_1}
&A_0\ar[l]_-{\p_0}
&Y_0=X\ar[l]_-{\i_0}
}
\leqno{\r:}
$$
satisfying
$$
\p_n\i_n=0
$$
in $\A$ for all $n=0,1,2,...$. The \emph{associated $X$-coaugmented cochain
complex} of such a sequence is then defined to be
$$
\xymat{
...
&A_n\ar[l]_-{\i_{n+1}\p_n}
&...\ar[l]_-{\i_n\p_{n-1}}
&A_1\ar[l]_-{\i_2\p_1}
&A_0\ar[l]_-{\i_1\p_0}
&X.\ar[l]_-{\i_0}
}
\leqno{C_\A(\r):}
$$

For an additive subcategory $\a\subset\A=\B\ho$, an $X$-coaugmented sequence
$\r$ as above will be called an \emph{$\a$-sequence} if $A_n$ belongs to
$\a$ for all $n$. Moreover it will be called \emph{$\a$-exact} if for any
object $A$ from $\a$, the induced sequence
$$
\Hom_\A(Y_{n+1},A)\to\Hom_\A(A_n,A)\to\Hom_\A(Y_n,A)
$$
is a short exact sequence of abelian groups for all $n\ge0$. Thus in this
case, the chain complex $C_\A(\r)$ is $\a$-exact in the sense of
\ref{1-complex}. In fact, for any object $A$ in $\a$ the differential
${\bar d}_n:\Hom_\A(A_{n+1},A)\to\Hom_\A(A_n,A)$ in $\Hom_\A(C_\A(\r),A)$
is then $\Hom_\A(\i_{n+1}\p_n,A)$, and one has
$$
\mathrm{ker}({\bar d}_{n-1})=\mathrm{im}({\bar d}_n)=\Hom_\A(Y_{n+1},A)
$$
for all $n$ and all $A\in\a$.
\end{Definition}

\begin{Proposition}\label{seqsec}
For each $X$-coaugmented sequence $\r$ in $\A$, any choice of
representatives $i_n\in\i_n$, $p_n\in\p_n$ in $\B_0$ and of tracks
$\alpha_n:p_ni_n\then0_{Y_n,Y_{n+1}}$ determines an $X$-coaugmented
secondary chain complex in $\B$ of the form
$$
\xymat{
\cdots
&A_3\ar[l]|-{i_4p_3}
&A_2\ar[l]|-{i_3p_2}\lluppertwocell<14>^{0}{^{\ \ \ \ \ \ \ i_4\alpha_3p_2}}
&A_1\ar[l]|-{i_2p_1}\lllowertwocell<-14>_{0}{_{i_3\alpha_2p_1\ \ \ \ \ \ \ }}
&A_0\ar[l]|-{i_1p_0}\lluppertwocell<14>^{0}{^{\ \ \ \ \ \ \ i_2\alpha_1p_0}}
&X.\ar[l]|-{i_0}\lllowertwocell<-14>_{0}{_{i_1\alpha_0\ \ \ \ \ \ }}
}
\leqno{C_\B(\r):}
$$
\end{Proposition}

\begin{proof}
Consider the diagram
$$
\includegraphics{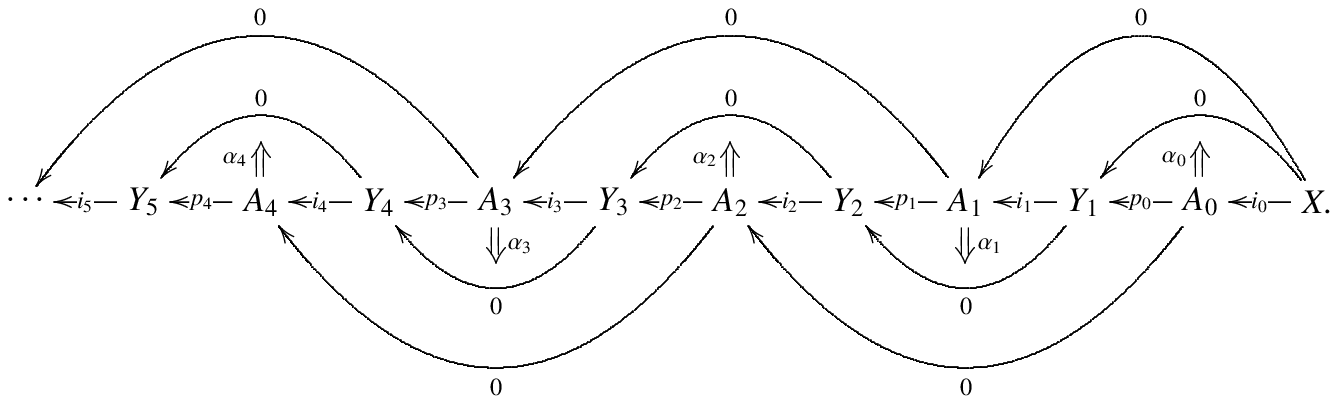}
$$
That this diagram yields on $C_\B(\r)$ above the structure of a secondary
chain complex, is equivalent to the identities
$$
i_{n+1}p_ni_n\alpha_{n-1}p_{n-2}=i_{n+1}\alpha_np_{n-1}i_{n-1}p_{n-2}.
$$
These are satisfied since one actually has
$$
p_ni_n\alpha_{n-1}=\alpha_np_{n-1}i_{n-1},
$$
as the next lemma shows.
\end{proof}

\begin{Lemma}
For any maps $f:X\to Y$, $f':Y\to Z$ and tracks $\alpha:f\then0_{X,Y}$,
$\alpha':f'\then0_{Y,Z}$ one has
$$
f'\alpha=\alpha'f.
$$
\end{Lemma}

\begin{proof}
This is a particular case of \ref{interchange}.
\end{proof}

\begin{Remark}
Strictly speaking, $C_\B(\r)$ depends on the choice of the $i_n$, $p_n$ and
$\alpha_n$; however it will be harmless in what follows to suppress these
from the notation.
\end{Remark}

\begin{Example}\label{tower}
Let $\B$ be a track category and suppose that $\A=\B\ho$ is equipped with
the structure of a triangulated category. Thus there is an endofunctor
$\RR:\A\to\A$ which is a self-equivalence, with an inverse equivalence
$\RR\1$, and one has a distinguished class of diagrams of the form
$$
A\ot B\ot C\ot\RR A,
$$
called exact triangles, which satisfy certain axioms. A \emph{fiber tower}
$\t$ over an object $X$ is a diagram in $\A$
\begin{equation}\label{fit}
\xymat{
X&X_0\ar@{=}[l]\ar[d]&X_1\ar[l]\ar[d]&X_2\ar[l]\ar[d]&...\ar[l]\\
&A_0&A_1&A_2&\cdots
}
\end{equation}
such that each $A_n\ot X_n\ot X_{n+1}\ot\RR A_n$ is an exact triangle in
$\A$. In particular, the composites $\RR\1X_{n+1}\ot A_n\ot X_n$ are zero
maps in $\A$. For an additive subcategory $\a$ in $\A$, call a fiber tower $\t$
\emph{$\a$-exact} if $A_i\in\a$ for all $i$ and moreover each of its exact
triangles induces a short exact sequence
$$
0\to\Hom(\RR\1X_{n+1},A)\to\Hom(A_n,A)\to\Hom(X_n,A)\to0
$$
for all $A\in\a$.

A fiber tower yields a system of coaugmented sequences in $\A$ of the form
\begin{equation}\label{seqseq}
\xymatrixrowsep{.1pc}
\begin{aligned}
\xymatrix{
   &         &                &                &A_0               &X\ar[l]\\
   &         &             A_1&       X_1\ar[l]&\RR A_0\ar[l]  &\RR X\ar[l]\\
A_2&X_2\ar[l]&\RR A_1\ar[l]&\RR X_1\ar[l]&\RR^2 A_0\ar[l]&\RR^2 X\ar[l]\\
   &         &                &\vdots
}
\end{aligned}
\end{equation}
which via delooping in $\A$ yields the $X$-coaugmented sequence
$$
\r(\t):\ \ ...\ot\RR^{-2}A_2\ot\RR^{-2}X_2\ot\RR\1A_1\ot\RR\1X_1\ot A_0\ot X.
$$
Thus by \ref{seqsec} each fiber tower over $X$ gives rise to an
$X$-coaugmented secondary chain complex $C_\B(\r(\t))$.
\end{Example}

\begin{Remark}\label{birgit}
Before the authors obtained the construction from \ref{seqsec}, a direct
topological proof that Adams resolutions give rise to a secondary complex
has been kindly provided to them by Birgit Richter \cite{Birgit}.
\end{Remark}

One then has

\begin{Theorem}\label{towerdif}
Assume either $\B$ is $\LL$-additive and $A\in\a$ implies $\LL A\in\a$, or $\B$ is
$\RR$-additive and $\RR$ preserves $\a$-exactness of complexes in
$\A$ (cf. \ref{hosec}). Then for any $\a$-exact fiber tower over an object $X$, any
$X$-coaugmented secondary chain complex associated to it (as in \ref{tower}
and \ref{seqsec})
is a $\b$-coresolution of $X$. Hence for any object $Y$ there is a
secondary differential
$$
d_{(2)}:\Ext^n_{\a\op}(X,Y)^m\to\Ext^{n+2}_{\a\op}(X,Y)^{m+1},
$$
where $\Ext^*_{\a\op}(X,Y)^m$ denotes either $\Ext^*_{\a\op}(\LL^mX,Y)$ or
$\Ext^*_{\a\op}(X,\RR^mY)$ in $\A\op$. The differential $d_{(2)}$ is well-defined by the
cohomology class $\brk{\B}\in H^3(\A;D)$ with $D$ in \ref{sigmaomega}.
\end{Theorem}

\begin{proof}
This follows directly from \ref{hosec}.
\end{proof}

\section{The $\mathrm{E}_3$ term of the Adams spectral sequence}

As in example \ref{staho} let $\A$ be the stable homotopy category of
spectra and let $\a\subset\A$ be the full subcategory of finite products of
Eilenberg-Mac Lane spectra over a fixed prime field $\F_p$. Let $X$ be a
spectrum of finite type, that is, for which the cohomology groups
$H^i(X;\F_p)$ are finite dimensional $\F_p$-vector spaces for all $i$. Then
the \emph{Adams fiber tower} of $X$ is given by
\begin{equation}\label{afit}
\xymat{
X&X_0\ar@{=}[l]\ar[d]&X_1\ar[l]\ar[d]&X_2\ar[l]\ar[d]&...\ar[l]\\
&H\wedge X_0&H\wedge X_1&H\wedge X_2&\cdots.
}
\end{equation}
Here $H=H\F_p$ is the Eilenberg-Mac Lane spectrum, the map $X_i\to H\wedge X_i$
is given by smashing $S^0\to H$ with $X_i$, and $X_{i+1}$ is the fiber of
this map. Since $X$ is of finite type all spectra $H\wedge X_i$ can be
considered to be objects of $\a$. By construction the Adams fiber tower is
$\a$-exact.

Since the category of spectra is a Quillen model category we know that $\A$
is the homotopy category of all spectra which are fibrant and cofibrant.
Using the cylinder of such spectra we obtain the additive track category
$\B$. That is, $\B$ consists of spectra which are fibrant and cofibrant,
of maps between such spectra, and tracks between such maps. Then $\B$ is
$\LL$-additive (and also $\Omega$-additive) and $A\in\a$ implies $\LL
A\in\a$. Therefore we can apply theorem \ref{towerdif} to the Adams fiber
tower where $\b$ is the full track subcategory of $\B$ on objects from
$\a$. Hence we get for a spectrum $Y$ the following diagram whose top row
is defined by any secondary $\b$-coresolution of $X$ and the bottom row is
the differential $d_{(2)}$ in the Adams spectral sequence.
\begin{equation}\label{d2comp}
\xymat{
\Ext^n_{\a\op}(X,Y)^m\ar[r]^{d_{(2)}}\ar[d]^\cong
&\Ext^{n+2}_{\a\op}(X,Y)^{m+1}\ar[d]^\cong\\
\Ext_{\mathscr A}^n(H^*X,H^*Y)^m\ar[r]^{d_{(2)}}
&\Ext_{\mathscr A}^{n+2}(H^*X,H^*Y)^{m+1}.
}
\end{equation}
The vertical isomorphisms in this diagram are defined in example
\ref{staho}.

\begin{Theorem}\label{coinc}
The diagram \ref{d2comp} commutes.
\end{Theorem}

This shows that $d_{(2)}d_{(2)}=0$ so that the secondary $\Ext$ in section
\ref{sext} coincides with the E$_3$ term of the Adams spectral sequence. In
the book \cite{Baues} a pair algebra $\mathscr B$ is computed which can be
used to describe algebraic models for secondary $\b$-coresolutions. This,
in fact, yields an algorithm computing the $d_{(2)}$ differential in the
Adams spectral sequence since we can use theorem \ref{coinc}.

\begin{proof}
In our terms the second differential of the Adams spectral sequence can be
understood in the following way: one is given a fiber tower $\t$ like
\ref{fit} or \ref{afit} over an object $X$ in the stable homotopy category, with the
associated $X$-coaugmented sequence $\r(\t)$ as in \ref{tower}. To it
corresponds by \ref{coaugdef} the associated $X$-coaugmented cochain
complex
$$
\begin{aligned}
\xymatrixcolsep{2pc}
\xymatrix{
...
&\Omega^{-n-1}A_{n+1}\ar[l]
&\Omega^{-n}A_n\ar[l]_-{\Omega^{-n}d^n}
&...\ar[l]
&\Omega\1A_1\ar[l]_-{\Omega\1d^1}
&A_0\ar[l]_-{d^0}
&X\ar[l]
}
\end{aligned}
\leqno{C_\A(\r(\t)):}
$$
where $d^n:A_n\to\Omega\1 A_{n+1}$ are the composites $A_n\to\Omega\1
X_{n+1}\to\Omega\1 A_{n+1}$ of maps in the exact triangles $X_{n+1}\to
X_n\to A_n\to\Omega\1X_{n+1}$ and
$A_{n+1}\to\Omega\1X_{n+2}\to\Omega\1X_{n+1}\to\Omega\1A_{n+1}$. Here all
$A_n$ are $\F_p$-module spectra, i.~e. Eilenberg-MacLane spectra of
$\F_p$-vector spaces, and moreover the sequences $X_n\to
A_n\to\Omega\1X_{n+1}$ are $\F_p$-exact, i.~e. applying $H^*(\_;\F_p)$ to
them yields short exact sequences. In particular, $H^*(C_\A(\r(\t));\F_p)$
is an $\mathscr A$-projective resolution of $H^*(X;\F_p)$.

Now choose new spectra $B_n$ fitting in exact triangles
$$
\xymat{
B_n\ar[r]^-{i^n}
&A_n\ar[r]^-{d^n}
&\Omega\1A_{n+1}\ar[r]
&\Omega\1 B_n
}
$$
and observe that by the octahedron axiom there is a commutative diagram of
(co)fibre sequences of the form
\begin{equation}\label{braid}
\begin{aligned}
\xymatrix@!C=2em{
\\
X_{n+2}\ar@/^5ex/[rr]\ar[dr]
&&X_n\ar@/^5ex/[rr]\ar[dr]
&&A_n\ar@/^5ex/[rr]\ar[dr]
&&\Omega\1A_{n+1}\ar@/^5ex/[rr]\ar[dr]
&&\cdots\\
&X_{n+1}\ar[ur]\ar[dr]
&&B_n\ar[ur]\ar[dr]
&&\Omega\1X_{n+1}\ar[ur]\ar[dr]
&&\Omega\1B_n\ar[ur]\ar[dr]
&\cdots\\
&&A_{n+1}\ar@/_5ex/[rr]\ar[ur]
&&\Omega\1X_{n+2}\ar@/_5ex/[rr]\ar[ur]
&&\Omega\1X_n\ar@/_5ex/[rr]\ar[ur]
&&\cdots\\
\ 
}
\end{aligned}
\end{equation}
so that in particular the original fiber tower $\t$ ``doubles'' to give two new
fiber towers $\t^{(2)}$ starting at $X_0$, resp. $X_1$, of the form
$$
\xymatrix{
...&X_n\ar[d]\ar[l]&X_{n+2}\ar[l]\ar[d]&X_{n+4}\ar[l]\ar[d]&...\ar[l]\\
\cdots&B_n&B_{n+2}&B_{n+4}&\cdots.
}
$$
The associated sequences $\r(\t^{(2)})$ and the cochain complexes
$$
\xymat{
...
&\Omega^{-n-2}B_{n+4}\ar[l]
&\Omega^{-n-1}B_{n+2}\ar[l]
&\Omega^{-n}B_n\ar[l]_-{\Omega^{-n}d(2)^n}
&...,\ar[l]
}
\leqno{C_\A(\r(\t^{(2)})):}
$$
where $d(2)^n:B_n\to\Omega\1B_{n+2}$ is the composite
$B_n\to\Omega\1X_{n+2}\to\Omega\1B_{n+2}$, are then obtained as in
\ref{tower}.

Let us now take any spectrum $Y$ and apply the stable homotopy classes
functor $\set{Y,\_}$ to the whole business. Because of the exact triangles
$B_n\to A_n\to\Omega\1A_{n+1}\to\Omega\1B_n$, there are isomorphisms
$$
\im\left(\set{Y,B_n}\to\set{Y,A_n}\right)\cong
\ker\left(\set{Y,A_n}\to\set{Y,\Omega\1A_{n+1}}\right).
$$
On the other hand it is known (see e.~g. \cite{Moss}) that the canonical maps
\begin{equation}\label{horep}
\set{Y,A_n}\xto\cong\Hom_{\mathscr A}(H^*(A_n;\F_p),H^*(Y;\F_p))
\end{equation}
are isomorphisms; it thus follows that the groups
$$
E_2^{s,t}(Y,X)=\frac
{\im\left(\set{Y,\Omega^{t-s}B_s}\xto{\set{Y,\Omega^{t-s}i_s}}\set{Y,\Omega^{t-s}A_s}\right)}
{\im\left(\set{Y,\Omega^{t-s+1}A_{s-1}}\xto{\set{Y,\Omega^{t-s+1}d^{s-1}}}\set{Y,\Omega^{t-s}A_s}\right)}
$$
are isomorphic to $\Ext_{\mathscr A}^s(H^*(X;\F_p),H^*(Y;\F_p))^t$.

Moreover (see again \cite{Moss}) the Adams differential $E_2^{s,t}\to
E_2^{s+2,t+1}$ is induced by the map
$$
\im\left(\set{Y,\Omega^{t-s}B_s}\to\set{Y,\Omega^{t-s}A_s}\right)
\to\im\left(\set{Y,\Omega^{t-s-1}B_{s+2}}\to\set{Y,\Omega^{t-s-1}A_{s+2}}\right)
$$
which sends the class of a stable map
$$
Y\to\Omega^{t-s}B_s\to\Omega^{t-s}A_s
$$
to the class of the composite
$$
Y\to\Omega^{t-s}B_s\xto{\Omega^{t-s}d(2)^s}\Omega^{t-s-1}B_{s+2}\to\Omega^{t-s-1}A_{s+2}
$$
or, which by \ref{braid} is the same, the composite
$$
Y\to\Omega^{t-s}B_s\to\Omega^{t-s-1}X_{s+2}\to\Omega^{t-s-1}A_{s+2}.
$$

To see then that the differential so defined coincides with the secondary
differential as constructed in \ref{indep}, \ref{sext} and \ref{towerdif},
let us choose zero tracks $\alpha_n$ for the composites $X_n\to
A_n\to\Omega\1X_{n+1}$ and switch from $C_\A(\r(\t))$ to the
$X$-coaugmented secondary cochain complex $C_\B(\r(\t))$ as defined in
\ref{seqsec}. Then according to \ref{indep}, given an element $\brk{c}$ of
$\Ext^s_{\mathscr A}(H^*(X;\F_p),H^*(Y;\F_p)))^t$, the corresponding
element $d_{(2)}\brk{c}\in\Ext^{s+2}_{\mathscr
A}(H^*(X;\F_p),H^*(Y;\F_p))^{t+1}$ is constructed in the following way.
First represent $\brk{c}$ by a cocycle in $C_\A(\r(\t))$, i.~e. by a
homomorphism of ${\mathscr A}$-modules $[c]:H^*(\Omega^{t-s}A_s;\F_p)\to
H^*(Y;\F_p)$ with $[c]\circ H^*(\Omega^{t-s}d^s;\F_p)=0$. By \ref{horep},
this homomorphism is in turn induced by a map $c:Y\to\Omega^{t-s}A_s$ such
that $d^s\circ c$ is nullhomotopic. Choosing a homotopy $\gamma:0\then
d^s\circ c$, according to \ref{indep} the class $d_{(2)}\brk{c}$ is
represented by the map
$Y\to\Omega^{t-s-1}A_{s+2}=\Omega\Omega^{t-s-2}A_{s+2}$ which corresponds
to the composite homotopy
$$
\xymatrix{
0=\Omega^{t-s-1}d^{s+1}\circ0\ar@{=>}[rr]^-{\Omega^{t-s-1}d^{s+1}\gamma}
&&\Omega^{t-s-1}d^{s+1}\circ\Omega^{t-s}d^s\circ c\ar@{=>}[r]^-{\delta c}
&0c=0
}
$$
from the zero map $Y\to\Omega^{t-s-2}A_{s+2}$ to itself, as in
$$
\xymat{
Y\ar[dd]_c\ddrtwocell<0>^0{_<4>\gamma}\\
\\
\Omega^{t-s}A_s\ar[r]_{\Omega^{t-s}d^s}\rrlowertwocell<-10.5>_{0}{_\hskip1ex\d}
&\Omega^{t-s-1}A_{s+1}\ar[r]^{\Omega^{t-s-1}d^{s+1}}
&\Omega^{t-s-2}A_{s+2}.
}
$$

Now according to the construction of $C_\B(\r(\t))$ given in \ref{seqsec},
this diagram reduces to the following diagram
\begin{equation}\label{ford2}
\begin{aligned}
\xymatrix@!C=5em{
Y\ar[dd]_c\ddrrtwocell<0>^0{_<3>\gamma}\\
\\
\Omega^{t-s}A_s\ar[r]
&\Omega^{t-s-1}X_{s+1}\ar[r]\rrlowertwocell<-10.5>_{0}{_\hskip1ex\alpha}
&\Omega^{t-s-1}A_{s+1}\ar[r]
&\Omega^{t-s-2}X_{s+2}\ar[r]
&\Omega^{t-s-2}A_{s+2}.
}
\end{aligned}
\end{equation}

Next note that because of the fibre sequence
$$
\Omega^{t-s}B_s\to\Omega^{t-s}A_s\to\Omega^{t-s-1}A_{s+1}\to\Omega^{t-s-1}B_s,
$$
choosing $\gamma:0\then\Omega^{t-s}d^s\circ c$ is equivalent to choosing a
lift of $c$ to a map $Y\to\Omega^{t-s}B_s$. Similar correspondences
between homotopies and liftings of maps take place further along the
sequence, as can be summarized in the following diagram
$$
\xymat{
Y\ar@{=}[d]\ar[r]
&\Omega^{t-s}B_s\ar[d]\ar[r]
&\Omega^{t-s-1}X_{s+2}\ar[d]\ar@{=}[r]
&\Omega^{t-s-1}X_{s+2}\ar[d]\ar[r]
&\Omega^{t-s-1}A_{s+2}\ar[d]\\
Y\ar[d]\ar[r]\drtwocell\omit{^\gamma}
&\Omega^{t-s}A_s\ar[d]\ar[r]\drtwocell\omit{=}
&\Omega^{t-s-1}X_{s+1}\ar[d]\ar[r]\drtwocell\omit{^\alpha}
&{*}\ar[d]\ar@{=}[r]\drtwocell\omit{=}
&{*}\ar[d]\\
{*}\ar[r]
&\Omega^{t-s-1}A_{s+1}\ar@{=}[r]
&\Omega^{t-s-1}A_{s+1}\ar[r]
&\Omega^{t-s-2}X_{s+2}\ar[r]
&\Omega^{t-s-2}A_{s+2},
}
$$
in which the columns form fiber sequences and the upper horizontal maps are
liftings corresponding to the homotopies indicated in lower squares. That
the resulting upper horizontal composite is indeed the lifting
corresponding to the composite homotopy in \ref{ford2} now follows from the
following standard homotopy-theoretic lemma which can be found e.~g. in
\cite{BauesAH}*{(2.9) on p. 263}:

\begin{Lemma}
Given a diagram
$$
\xymat{
F\ar[d]\ar@{-->}[r]
&F'\ar[d]\ar@{-->}[r]
&F''\ar[d]
\\
E\ar[d]\ar[r]^e\drtwocell\omit{^\gamma}
&E'\ar[d]\ar[r]\drtwocell\omit{^\delta}
&E''\ar[d]
\\
B\ar[r]
&B'\ar[r]^b
&B''}
$$
whose columns are fiber sequences and upper horizontal maps are liftings
corresponding to the indicated homotopies, then the composite $F\to F''$ is
the lifting corresponding to the composite homotopy $\delta e\comp
b\gamma$.
\end{Lemma}

\end{proof}

\section*{Acknowledgements}

The authors are grateful to Birgit Richter for \ref{birgit} and to Gerald
Gaudens for discussions on several aspects of the paper.

\end{document}